\documentclass[11pt]{article}

\usepackage{amsfonts,amsthm}
\usepackage{amsmath,mathtools}   
\usepackage{latexsym}
\usepackage{amssymb}
\usepackage{mathrsfs}
\usepackage{pifont}
\usepackage{tikz}
\usepackage[all]{xy}
\usepackage[hidelinks]{hyperref}

\newtheoremstyle{newrem}{3pt}{3pt}{}{}
{\bfseries}{.}{.5em}{}

\newtheorem{theo}{Theorem}[section]
\newtheorem{thm}{Theorem}[section]

\newtheorem*{theo*}{Theorem}

\newtheorem{prop}[theo]{Proposition}

\newtheorem{coro}[theo]{Corollary}

\theoremstyle{newrem}
\newtheorem{remax}[theo]{Remark}

\theoremstyle{definition}
\newtheorem{defi}[theo]{Definition}
\newtheorem{exam}[theo]{Example}

\newtheorem*{term*}{Notation/Terminology}

\newcommand{\bC}{\mathbb{C}}
\newcommand{\bZ}{\mathbb{Z}}

\newcommand{\tm}{\textbf{m}}

\newcommand{\cA}{\mathcal{A}}
\newcommand{\cB}{\mathcal{B}}

\newcommand{\oT}{\overline{T}}

\usepackage{geometry}
\begin{document}

\title{The missing label of $\mathfrak{su}_3$ and its symmetry}
\date{}

\author{Nicolas Cramp\'e\footnote{Institut Denis-Poisson CNRS/UMR 7013 - Universit\'e de Tours - Universit\'e d'Orl\'eans, Parc de Grandmont, 37200 Tours, France. crampe1977@gmail.com }
, 
Lo\"ic Poulain d'Andecy\footnote{ Laboratoire de math\'ematiques de Reims UMR 9008, Universit\'e de Reims Champagne-Ardenne,
Moulin de la Housse BP 1039, 51100 Reims, France. loic.poulain-dandecy@univ-reims.fr}
,
Luc Vinet\footnote{Centre de recherches math\'ematiques, Universit\'e de Montr\'eal, P.O. Box 6128, Centre-ville Station, Montr\'eal (Qu\'ebec), H3C 3J7, Canada. vinet@CRM.UMontreal.ca}
\footnote{Insitut de valorisation des donn\'ees (IVADO), Montr\'eal (Qu\'ebec), H2S 3H1, Canada. }}

\maketitle

\begin{abstract}
We present explicit formulas for the operators providing missing labels for the tensor product of two irreducible representations of $\mathfrak{su}_3$. The result is seen as a particular representation of the diagonal centraliser of $\mathfrak{su}_3$ through a pair of tridiagonal matrices. Using these explicit formulas, we investigate the symmetry of this missing label problem and we find a symmetry group of order 144 larger than what can be expected from the natural symmetries. Several realisations of this symmetry group are given, including an interpretation as a subgroup of the Weyl group of type $E_6$, which appeared in an earlier work as the symmetry group of the diagonal centraliser. Using the combinatorics of the root system of type $E_6$, we provide a family of representations of the diagonal centraliser by infinite tridiagonal matrices, from which all the finite-dimensional representations affording the missing label can be extracted.
Besides, some connections with the Hahn algebra, Heun--Hahn operators and Bethe ansatz are discussed along with some similarities with the well-known symmetries of the Clebsch--Gordan coefficients.
\end{abstract}

\vspace{3mm}

\section*{Introduction}

We consider two irreducible finite-dimensional representations $[m_1,m_2]$ and $[m'_1,m'_2]$ of the Lie algebra $\mathfrak{su}_3$, corresponding to two highest weights. 
We form the tensor product $[m_1,m_2]\otimes [m'_1,m'_2]$ and we consider a third representation $[m''_1,m''_2]$ appearing in the Clebsch--Gordan series:
\begin{equation}\label{CG-series}
[m_1,m_2]\otimes [m'_1,m'_2]=\bigoplus_{(m''_1,m''_2)} [m''_1,m''_2]^{\oplus d_{\tm}}\ .
\end{equation}
We denote collectively by $\tm=(m_1,m_2,m'_1,m'_2,m''_1,m''_2)$ the three highest weights. The multiplicity $d_{\tm}$ with which the representation $[m''_1,m''_2]$ appears in this decomposition is the Littlewood--Richardson coefficient. This multiplicity may be larger than one and this is where a missing label problem appears. Namely, the natural labels can specify a vector in the representation $[m''_1,m''_2]$ but cannot distinguish between the different copies of this representation appearing in the series.

To resolve this ambiguity, the centraliser in $U(\mathfrak{su}_3)\otimes U(\mathfrak{su}_3)$ of the diagonal embedding of $U(\mathfrak{su}_3)$ must be considered.
In \cite{CPV1}, two generators $X$ and $Y$ of this algebra have been identified and an algebraic description of $Z_{2}(\mathfrak{su}_3)_{\tm}$, 
depending on the 6 parameters in $\tm$ is obtained.
The algebra $Z_{2}(\mathfrak{su}_3)_{\tm}$ acts on the different copies of $[m''_1,m''_2]$ in (\ref{CG-series}) and is the relevant algebra providing missing label operators. 
In fact, for easy future reference, we will call \emph{physical} this representation of $Z_{2}(\mathfrak{su}_3)_{\tm}$ since this algebra admits other  
representations not associated with the missing label problem.

The action of $X$ in the physical representation was calculated in \cite{PST}. It follows from this calculation that $X$ can distinguish between the various copies of $[m''_1,m''_2]$ and 
thus its eigenvalues can indeed serve as a missing label. The second operator $Y$, which allows (through various combinations of $X$ and $Y$) to select among different missing label operators, 
was not much considered so far. Moreover, the matrix obtained in \cite{PST} for $X$ did not display symmetries apart from the obvious ones. 
The first main result of this paper provides the action of $X$ and $Y$ in a physical representation (in a different basis than the one in \cite{PST}) 
in such a way that symmetries will appear (see Section \ref{sec-XY} and Appendix \ref{app:XY}).
By a symmetry of the missing label, we mean a transformation on the parameters $\tm$ such that the matrices for $X$ and $Y$ are transformed 
to equivalent matrices (up to a $\pm$ sign for $X$). Thus the eigenvalues will be preserved (up to a sign for $X$) under these transformations. 

We prove the following theorem.
\begin{theo*}
On the parameters in $\tm$, the symmetries of the missing label are described as transformations of the following arrangement:
\begin{equation}
 \label{eq:thm1}
\left[\begin{array}{ccc|ccc}
m_1 & m'_1 & m''_2                             \ \ \ & \ \ \ m_1+\ell-n & m'_1+\ell-n & m''_2+\ell-n  \\[0.5em]
m'_1+m'_2-\ell & m_1+m_2-\ell & n      \ \ \  & \ \ \ m'_1+m'_2-n & m_1+m_2-n & \ell \\[0.5em]
m_2+n-\ell & m'_2+n-\ell & m''_1+n-\ell \ \ \ & \ \ \ m_2 & m'_2 & m''_1
\end{array}\right]
\end{equation}
which are combinations of the following operations:\\
$\bullet$ simultaneously on the left and right $3\times 3$ squares: permutations of the lines, permutations of the columns and transposition of matrices;\\
$\bullet$ exchange of the left and right $3\times 3$ squares.\\
The following combinations $\ell$ and $n$ of the parameters $\tm$ has been used in \eqref{eq:thm1} and will appear everywhere throughout the paper:
\begin{equation}\label{eq:nl2}
 \ell=\frac{1}{3}(m_1+2m_2+m'_1+2m'_2-m''_1-2m''_2)\ \ \ \ \text{and}\ \ \ \ n=\frac{1}{3}(2m_1+m_2+2m'_1+m'_2-2m''_1-m''_2)\ .
\end{equation}
They satisfy also the formulas:
\begin{equation}\label{eq:nl1}
 m''_1=m_1+m'_1+\ell-2n\ \ \ \ \ \text{and}\ \ \ \ \ \ m''_2=m_2+m'_2+n-2\ell\ .
\end{equation}
\end{theo*}
\noindent To be precise, the transformations in the theorem transform $X$ to a matrix equivalent to $\pm X$, the sign being determined by the number of transpositions of lines and of columns, and they transform $Y$  to a matrix equivalent to $Y$. Note that the two $3\times 3$ squares above are magic squares, in the sense that the sum over any line or column or diagonal or antidiagonal gives the same result. 
The exchange of the left and right $3\times 3$ squares commutes with all other symmetries, so the whole transformation group is the direct product of a group of order 72 
with $\bZ/2\bZ$. The subgroup of order 72 acting simultaneously on both $3\times 3$ squares has the same structure as the 
group describing the symmetries of the Clebsch--Gordan coefficients, 
or $3j$-symbols, for $\mathfrak{su}_2$. We describe the similarities between both 
groups of symmetries in Section \ref{sec:CG}. An alternative ways to describe these symmetries is also given with affine transformations on a torus (see Figure \ref{fig}).

The Littlewood--Richardson coefficient $d_{\tm}$ is the minimum of the $18$ numbers appearing in the array above. Therefore, the condition (selection rules) for $\tm$ to represent a 
physical situation in the Clebsch--Gordan series is that every number in this array be a positive integer. Of course a symmetry of the missing label preserves 
in particular the size of the matrices $X$ and $Y$ (namely, the dimension of the physical representation), 
and thus is in particular a symmetry of the Littlewood--Richardson coefficient $d_{\tm}$. Indeed it is clear that the transformations of the theorem are symmetries of 
the Littlewood--Richardson coefficient (they permute the 18 elements). We recover the combinatorial results of \cite{BR}.

In addition to the symmetries of the matrices $X$ and $Y$, 
two other results giving insights in the description of the missing labels of $\mathfrak{su}_3$ are provided. In Section \ref{sec:CGsu2}, a connection with the algebraic 
Hahn--Heun operator is proven: the matrices $X$ is a combination of the Hahn bispectral operators of the eponymous orthogonal polynomials which arise in the study of the 
Clebsch--Gordan coefficients for $\mathfrak{su}_2$. Then, in Section \ref{sec:gm}, a link with integrable systems is pointed out and Bethe ansatz techniques are used to 
give the eigenvalues of $X$ in terms of the Bethe roots. Different examples of eigenvalues of $X$ are given by a direct computation and by the Bethe ansatz to 
prove conjectures given in \cite{CaSt}.

To give an other way of understanding the symmetries of $X$ and $Y$, we review the algebraic 
description of $Z_{2}(\mathfrak{su}_3)_{\tm}$ in Section \ref{sec-alg} as well as its symmetries \cite{CPV1}.
Recall that the parameters $\tm$ correspond to three weights of $\mathfrak{su}_3$. There is a natural embedding of the root system of $\mathfrak{su}_3^{\oplus 3}$ into the root system of type $E_6$. Through this embedding, there is a natural action of the group $\pm W(E_6)$ on the set of weights of $\mathfrak{su}_3^{\oplus 3}$ (where $W(E_6)$ is the Weyl group and the $\pm$ sign accounts 
for the additional central symmetry $-\text{Id}$). So this group $\pm W(E_6)$ acts on our parameters $\tm$, and remarkably, it was shown in \cite{CPV1} that we obtain this way a group of symmetries of the algebra $Z_{2}(\mathfrak{su}_3)_{\tm}$ (see Theorem \ref{th:symE6}).
The symmetries of the missing label are identified within this larger transformation group in Section \ref{subsec-symE6}.

Finally, we describe in Section \ref{sec:rep} some representations of the missing label algebra $Z_{2}(\mathfrak{su}_3)_{\tm}$. Remarkably, there exists a geometric way to obtain a 
set of representations of this algebra by using the polytope associated to the root system $E_6$. More precisely, to some $4$-faces of this polytope, we associate an infinite-dimensional representation of $Z_{2}(\mathfrak{su}_3)_{\tm}$, where $X$ and $Y$ are infinite tridiagonal matrices. Finite-dimensional representations can be extracted from the infinite-dimensional ones and the physical representations are identified among this set of representations.

\section{Formulas for the missing label operators}\label{sec-XY}

\subsection{The Lie algebra $\mathfrak{su}_3$ and its irreducible representations}

The complexified Lie algebra $\mathfrak{su}_3$ is generated by the elements $\{\, e_{pq}\, |\, 1\leq p,q\leq 3\, \}$ which satisfy
the following defining relations, for $1\leq p,q,r,s\leq 3$,
\begin{eqnarray}\label{eq:defglN}
&& [e_{pq}, e_{rs} ] = \delta_{qr} e_{ps} - \delta_{sp} e_{rq}\quad \text{and}\qquad \sum_{i=1}^3 e_{i i}=0\,.
\end{eqnarray}
We define also the Cartan generators by $h_1=e_{11}-e_{22}$ and  $h_2=e_{22}-e_{33}$.

The finite-dimensional irreducible representations $[m_1,m_2]$ of $\mathfrak{su}_3$ correspond to a dominant integral weight
(see \textit{e.g.} \cite{FH} for the representation theory of $\mathfrak{su}_3$).
This consists of two integers $m_1,m_2 \in \mathbb{Z}_{> 0}$ such that the highest weight vector $v_{m_1,m_2}$ satisfies:
\begin{eqnarray}
 h_{p}v_{m_1,m_2}&=& (m_p-1) v_{m_1,m_2} \quad \text{with} \quad p=1,2\label{weight-vector}\\
 e_{pq}v_{m_1,m_2}&=&0 \quad \text{with} \quad 1\leq p< q \leq 3 \ .\label{weight-vector2}
\end{eqnarray}
With an abuse of terminology, we will speak of a vector of weight $(m_1,m_2)$ for a vector $v_{m_1,m_2}$ satisfying \eqref{weight-vector}-\eqref{weight-vector2}. 
The representation $[m_1,m_2]$ is usually represented by a Young diagram with $m_1+m_2-2$ boxes in the first row and $m_2-1$ boxes in the second row.

The decomposition of the tensor product of two dominant integral weights $[m_1,m_2]$ and $[m'_1,m'_2]$ is given by:
\begin{equation}
[m_1,m_2] \otimes [m_1',m_2'] = \bigoplus_{m''_1,m''_2\in \mathbb{Z}_{> 0}}  \  [m''_1,m''_2]^{\oplus d_{\tm}} \ ,  \label{eq:ds}
\end{equation}
where the integer $d_{\tm}\geq 0 $ is the multiplicity, 
called the Littlewood--Richardson coefficient and $\tm=(m_1,m_2,m'_1,m'_2,m''_1,m''_2)$.
The component in this decomposition corresponding to the weight $ [m''_1,m''_2]$ can be written as a vector space:
\begin{equation}
 [m''_1,m''_2]^{\oplus d_{\tm}} =M_{\tm} \otimes [m''_1,m''_2]\,,
\end{equation}
where $M_{\tm}$ is the multiplicity space. This space $M_{\tm}$ can be identified with the subspace of 
$[m_1,m_2] \otimes [m_1',m_2']$ consisting of the highest-weight vectors of weight $(m''_1,m''_2)$.

The dimension of the multiplicity space $M_{\tm}$ is the multiplicity $d_{\tm}$. If non-zero, the multiplicity $d_{\tm}$ is the minimum of the following set (see \cite{PST} or \cite{OR}):
\begin{equation}\label{eq:deg}\left\{\begin{array}{cccccc}
m_1\,, & m'_1\,, & m''_2\,,                             & m_1+\ell-n\,, & m'_1+\ell-n\,, & m''_2+\ell-n\,,\!\!  \\[0.5em]
\!\!m'_1+m'_2-\ell\,,\!\! & m_1+m_2-\ell\,,\!\! & n\,,      & m'_1+m'_2-n\,,\!\! & m_1+m_2-n\,,\!\! & \ell\,, \\[0.5em]
m_2+n-\ell\,, & m'_2+n-\ell\,, & m''_1+n-\ell\,,\!\! & m_2\,, & m'_2\,, & m''_1
\end{array}\right\}\,,
\end{equation} 
where $n$ and $\ell$ are defined by \eqref{eq:nl2} or \eqref{eq:nl1}. 
The $18$ numbers are arranged as in the two magic squares in the Introduction. 
Note that the $9$ elements on the right are the $9$ elements on the left plus $\ell-n$. So if $n\leq \ell$ the multiplicity is the minimum of the left $9$ elements, 
while if $\ell\leq n$ the multiplicity is the minimum of the right $9$ elements. If any of the $18$ elements above is negative or zero, 
then the multiplicity is $0$ and we have no apparition of $[m''_1,m''_2]$ in $[m_1,m_2] \otimes [m_1',m_2']$.

\subsection{Polarised traces and the diagonal centraliser\label{sec:ptt}} 

\paragraph{Casimir elements.} The Casimir elements of $U(\mathfrak{su}_3)$, 
the universal enveloping algebra of $\mathfrak{su}_3$, are:
\begin{equation}\label{eq:KL}
K= e_{i_1 i_2} e_{i_2 i_1}\ \ \ \ \text{and}\ \ \ \ L=  \frac{1}{2}(e_{i_1 i_2} e_{i_2 i_3}e_{i_3i_1}+e_{i_2 i_1} e_{i_3 i_2}e_{i_1i_3})\ ,
\end{equation}
where the sums from $1$ to $3$ over the repeated indices are understood. These two elements generate the center of $U(\mathfrak{su}_3)$. 

In the finite-dimensional irreducible representation $[m_1,m_2]$ of $\mathfrak{su}_3$, the values of these Casimir elements are respectively
\begin{equation}\label{value-Cas}
k_{m_1,m_2}=\frac{2}{3}( m_1^2+ m_2^2+ m_1 m_2)-2\ \ \ \ \ \text{and}\ \ \ \ \ l_{m_1,m_2}=\frac{1}{9}(m_1+2m_2)(2m_1+m_2)(m_1-m_2)\ .
\end{equation}
The choice of $K$ and $L$ was dictated by symmetry. Indeed, the automorphism of $U(\mathfrak{su}_3)$ defined by
\begin{equation}\label{sigma}
\sigma\ :\ e_{ij} \mapsto - e_{\overline{\jmath}\,\overline{\imath}}\,,\qquad \text{where} \quad\overline{1}=3\,,\ \overline{2}=2\,,\ \overline{3}=1\,,
\end{equation}
transforms the Casimir elements as follows:
\begin{equation}
 \sigma(K)=K\ \ \ \ \ \text{and}\ \ \ \ \ \sigma(L)=-L\ .
\end{equation}
For representations, the automorphism $\sigma$ is related to the operation of taking the dual, that is, it replaces the highest-weight $(m_1,m_2)$ by $(m_2,m_1)$.

\paragraph{Polarised traces.}
Let us introduce the following elements of $U(\mathfrak{su}_3)^{\otimes 2}$, for $1\leq p,q\leq 3$:
\begin{equation}\label{eijeps}
e_{pq}^{(1)}= e_{pq} \otimes 1 \,, \quad 
e_{pq}^{(2)}=1 \otimes e_{pq}\ \quad\text{and}\quad\ \delta(e_{pq})=e_{pq}^{(1)}+e_{pq}^{(2)}\ .
\end{equation}
The map $\delta$ from $U(\mathfrak{su}_3)$ to $U(\mathfrak{su}_3)^{\otimes 2}$  extends to a morphism of algebra, called the diagonal embedding.

The polarised traces are the following elements of $U(\mathfrak{su}_3)^{\otimes 2}$, for $ a_1,\dots,a_\mu=1,2$,
\begin{equation}
 T^{(a_1,\dots,a_\mu)}=e_{i_2i_1}^{(a_1)}e_{i_3i_2}^{(a_2)}\dots e_{i_1i_\mu}^{(a_\mu)}\ \ \ \ \ \ \ \text{and}\ \ \ \ \ \ \ 
 \oT^{(a_1,\dots,a_\mu)}=e_{i_1i_2}^{(a_1)}e_{i_2i_3}^{(a_2)}\dots e_{i_\mu i_1}^{(a_\mu)}
\end{equation}
where again the sum signs are omitted. Let us remark that $T^{(1,1)}=\oT^{(1,1)}$ and $\oT^{(1,1,1)}=T^{(1,1,1)}+3T^{(1,1)}$.

\paragraph{Generators of the diagonal centraliser.} Let $Z_2(\mathfrak{su}_3)$ denote the centraliser in $U(\mathfrak{su}_3)^{\otimes 2}$ of the diagonal embedding 
of $U(\mathfrak{su}_3)$, that is, the subalgebra of all elements in $U(\mathfrak{su}_3)^{\otimes 2}$ which commute with $\delta(e_{pq})=e_{pq}^{(1)}+e_{pq}^{(2)}$ for $1\leq p,q \leq 3$.
In \cite{CPV1}, it is shown, using classical invariant theory, that the following elements
\begin{equation}\label{defgen}
\begin{array}{l}
k_1=K\otimes 1\,,\ \ k_2=1\otimes K\,,\ \ k_3=\delta(K)\,,\ \  l_1=L\otimes 1\,,\ \ l_2=1\otimes L\,,\ \ l_3=\delta(L)\,,\\[0.5em]
\displaystyle X=\frac{1}{2}(T^{(1,1,2)}-T^{(1,2,2)})+\frac{1}{3}(l_1-l_2)\,,\\[0.8em]
\displaystyle Y=\frac{1}{2}(T^{(1,1,2,2)}+\oT^{(1,1,2,2)})-\frac{1}{12}(T^{(1,2)})^2-\frac{5}{12}T^{(1,1)}T^{(2,2)}-2T^{(1,2)}\,,
\end{array}\end{equation}
generate the whole centraliser $Z_2(\mathfrak{su}_3)$. In terms of the polarised traces, we have $k_i=T^{(i,i)}$ and $l_i=T^{(i,i,i)}+\frac{3}{2}k_i$ for $i=1,2$, and $k_3=k_1+k_2+T^{(1,2)}$ and $l_3=l_1+l_2+3(T^{(1,1,2)}+T^{(1,2,2)})+9T^{(1,2)}$.

The centraliser acts on the multiplicity space $M_{\tm}$, 
which is thus a representation of $Z_2(\mathfrak{su}_3)$. In particular, $X$ and $Y$ have a well-defined representation on the multiplicity space $M_{\tm}$, denoted 
by $X_{\tm}$ and $Y_{\tm}$ and called the \emph{physical representation} associated to $\tm$.

\begin{remax}
These elements of the centraliser were chosen in \cite{CPV1} so that they behave nicely under the symmetry of the algebra $Z_2(\mathfrak{su}_3)$ 
(see Section \ref{sec-alg} and Appendix \ref{app:XY}). In particular, the exchange of the left and right copies of $U(\mathfrak{su}_3)$ in $U(\mathfrak{su}_3)^{\otimes 2}$ induces the following transformation of the generators:
\[k_1\leftrightarrow k_2\,,\ \ k_3\mapsto k_3\,,\ \ \ \ l_1\leftrightarrow l_2\,,\ \ l_3\mapsto l_3\,,\ \ \ X\mapsto -X\,,\ \ \ Y\mapsto Y\ ,\]
while the automorphism $\sigma\otimes\sigma$, where $\sigma$ is defined in (\ref{sigma}), acts as follows
\[k_1,k_2,k_3\mapsto k_1,k_2,k_3\,,\ \ \ \ l_1,l_2,l_3\mapsto -l_1,-l_2,-l_3\,,\ \ \ \ \ X\mapsto -X\,,\ \ \ Y\mapsto Y\ .\]
This is easily deduced from the fact that $\sigma\otimes\sigma$ sends the polarised trace $T^{(a_1,\dots,a_\mu)}$ to $(-1)^\mu \oT^{(a_1,\dots,a_\mu)}$ 
and the formulas $\oT^{(1,1,2)}=T^{(1,1,2)}+3 T^{(1,2)}$ and $\oT^{(1,2,2)}=T^{(1,2,2)}+3 T^{(1,2)}$.
\end{remax}

\begin{remax}\label{rem-Cas}
In the representation $M_{\tm}$ of $Z_2(\mathfrak{su}_3)$, the Casimir elements $k_1,\, k_2,\, k_3,\ l_1,\,l_2,\,l_3,$ act as numbers, see (\ref{value-Cas}). 
The couple of numbers for $(k_1,l_1)$ (resp. $(k_2,l_2)$ and $(k_3,l_3)$)  is $(k_{m_1,m_2},l_{m_1,m_2})$ (resp.  $(k_{m'_1,m'_2},l_{m'_1,m'_2})$ and $(k_{m''_1,m''_2},l_{m''_1,m''_2})$).
So the multiplicity space can be seen as a representation of a specialisation of $Z_2(\mathfrak{su}_3)$, denoted $Z_2(\mathfrak{su}_3)_{\tm}$, 
where the Casimir elements are replaced by their values. This algebra is studied in details in Section \ref{sec-alg}.
\end{remax}

\subsection{Explicit formulas for $X_{\tm}$ and $Y_{\tm}$ \label{subsec-formulas}}

In \cite{PST}, a formula for the action of $X_{\tm}$ on the multiplicity space $M_{\tm}$ is already given (the operator $Y_{\tm}$ is not considered). 
In the basis chosen in \cite{PST}, $X_{\tm}$ takes the form of a quasi-tridiagonal matrix 
and its symmetry properties presented in the next section are not manifest. The same would be observed for the action $Y_{\tm}$ in this basis if it was somehow calculated.
It turns out that there exists a basis of the multiplicity space $M_{\tm}$ where $X_{\tm}$ and $Y_{\tm}$ are simultaneously tridiagonal matrices.  
The explicit form of the vectors of this basis is given in the Appendix \ref{app:XY} together with details regarding the construction of the representations of $X_{\tm}$ or $Y_{\tm}$ given below.

Let us define
\[(\xi_1\,,\ \xi_2\,,\ \xi_3)=\left\{\begin{array}{l}
(\ell\,,\ m_2\,,\ m'_1+\ell-n)\ \ \ \text{if $\ell\leq n$\,,}\\[0.5em]
(n\,,\ m'_1\,,\ m_2+n-\ell)\ \ \ \text{if $n\leq \ell$\,,}
\end{array}\right.\ \ \ \ \text{and}\ \ \ \ (\xi_4\,,\ \xi_5\,,\ \xi_6)=(\ell-m'_2\,,\ n-m_1\,,\ 0)\,,\]
as well as the following linear combinations of these parameters
\[\Lambda=\sum_{i=1}^3\xi_i-\sum_{i=4}^6\xi_i+2|\ell-n|\ \ \ \ \text{and}\ \ \ \lambda_{\pm}=\pm\frac{1}{2}\Lambda+\frac{1}{6}\sum_{i=1}^6\xi_i\ .\]
With these notations, the multiplicity $d_{\tm}$ given by \eqref{eq:deg}, assuming that $d_{\tm}>0$, can be rewritten as follows 
\[d_{\tm}=\xi_a-\xi_b\,,\ \ \ \ \ \quad\text{where\ \ \ $\xi_a=\text{min}\{\xi_1,\xi_2,\xi_3\}\ $ and $\ \xi_b=\text{max}\{\xi_4,\xi_5,\xi_6\}$\ .}\]
The matrices $X_{\tm}$ and $Y_{\tm}$ are given by the following tridiagonal matrices:
\begin{equation}\label{eq:XY}
 X_{\tm}=\left(\begin{array}{ccccc}
    a_{11} & a_{12}  &0  & \dots  & 0  \\
    a_{21} & a_{22}  & \ddots  & \ddots & \vdots  \\
     0 & \ddots  & \ddots & \ddots & 0 \\
     \vdots & \ddots & \ddots & \ddots & a_{d-1,d}\\
     0 & \dots & 0 &    a_{d,d-1}  &  a_{d,d}
 \end{array}\right)\ \ \ \ \text{and}\ \ \ \ Y_{\tm}=\left(\begin{array}{ccccc}
    b_{11} & b_{12}  &0  & \dots  & 0  \\
    b_{21} & b_{22}  & \ddots  & \ddots & \vdots  \\
     0 & \ddots  & \ddots & \ddots & 0 \\
     \vdots & \ddots & \ddots & \ddots & b_{d-1,d}\\
     0 & \dots & 0 &    b_{d,d-1}  &  b_{d,d}
 \end{array}\right)
\end{equation}
where $d=d_{\tm}$ and the off-diagonal coefficients are:
\[\begin{array}{rcl}a_{j,j+1}=(j+\xi_b-\xi_1)(j+\xi_b-\xi_2)(j+\xi_b-\xi_3) \quad & \text{and} & \quad b_{j,j+1}=a_{j,j+1}(j+\xi_b-\lambda_-)\ ,\\[0.5em]
a_{j+1,j}=(j+\xi_b-\xi_4)(j+\xi_b-\xi_5)(j+\xi_b-\xi_6) \quad & \text{and} & \quad b_{j+1,j}=a_{j+1,j}(j+\xi_b-\lambda_+)\ .
\end{array}\]
To write the diagonal coefficients, we set:
\[E_{j,k}=\sum_{i=1}^6x_i(j)^k\,,\ \ \ \ \text{where $x_i(j)=(\xi_i-j-\xi_b+\frac{1}{2})$\ for $i\in\{1,\dots,6\}$\ .}\]
Then we have:
\[a_{jj}=-\frac{1}{108}\Bigl(\frac{7}{2}E_{j,1}^3-18 E_{j,1}E_{j,2}+18 E_{j,3}\Bigr)-\frac{1}{24}E_{j,1}\Bigl(\Lambda^2+2\Bigr)\ ,\]
\[b_{jj}=\frac{1}{288}\Bigl(\frac{5}{2}E_{j,1}^4+32 E_{j,1}E_{j,3}+6\bigl(E_{j,2}^2-3E_{j,1}^2E_{j,2}-4E_{j,4})+6E_{j,2}(\Lambda^2+2)-3E_{j,1}^2(\Lambda^2-2)-\frac{3}{2}\Lambda^4+6\Lambda^2-36\Bigr)\ .\]
The off-diagonal coefficients of $X_{\tm}$ and $Y_{\tm}$ take a very simple form whereas their diagonal coefficients are more involved.
The choice of their explicit expressions is motivated by the fact that some symmetries are easily extracted as will be seen in the next section.

\begin{remax}
Using the symmetries that shall be described in Section \ref{sec-sym}, it is observed that, for example, 
$\ell$ can be always chosen as the minimum of the 18 elements in (\ref{eq:deg}). In this case, we have $d=\ell$, $(\xi_1,\xi_2,\xi_3)=(\ell,\,m_2,\,m'_1+\ell-n)$ and $\xi_b=0$.
\end{remax}

It is easy to check that the non-diagonal coefficients of $X_{\tm}$ in (\ref{eq:XY}) are non-zero. This implies that the action of $X_{\tm}$ on $M_{\tm}$ is cyclic 
(the whole space $M_{\tm}$ is obtained from repeated action of $X_{\tm}$ on the first basis vector for example) and that $X_{\tm}$ is diagonalisable (since it can be symmetrized). 
As explained in \cite{PST} using a different formula for $X_{\tm}$, this allows to conclude that the eigenvalues of $X_{\tm}$ on $M_{\tm}$ are all distinct. 
Therefore they can serve as labels for the different copies of $[m''_1,m''_2]$ appearing in the tensor product $[m_1,m_2] \otimes [m_1',m_2']$.
In this sense the operator $X_{\tm}$ provides a solution of the ``missing label'' problem. 

A direct investigation shows that $\lambda_-\leq \xi_b$ or $\lambda_{-}> \xi_a$ (and similarly for $\lambda_+$) which implies that $b_{j,j+1}b_{j+1,j}$ are non-zero for $1\leq i \leq d_\tm-1$.
Then, the same reasoning as above shows that the eigenvalues of $Y_{\tm}$ on $M_{\tm}$ are all distinct and can be also used as a missing label.

In summary, the values of the Casimir operators $k_3$ and $l_3$ of the 
diagonal embedding of $\mathfrak{su}_3$ determine $m''_1$ and $m''_2$. 
The eigenvalues of the operator $X_{\tm}$ (or $Y_{\tm}$) identifies the different copies of $[m''_1,m''_2]$ in the decomposition of $[m_1,m_2] \otimes [m_1',m_2']$.
Then, the diagonal images of the labelling operators of vectors inside an irreducible representation 
of $\mathfrak{su}_3$ (\textit{e.g.} the hypercharge, the isospin and a projection of the isospin) determine uniquely the vector in  $[m''_1,m''_2]$. 

We note that there does not seem to exist exact formulas, in the general case, for the eigenvalues of $X_{\tm}$ or $Y_{\tm}$. 
In fact, diagonalising $X_{\tm}$ in some examples immediately puts $Y_{\tm}$ in a complicated form (and vice versa). So it seems that the above 
simultaneous tridiagonal form for $X_{\tm}$ or $Y_{\tm}$ is the simplest presentation that can be achieved, and one may find it quite remarkable that such a simple form exists.
In Section \ref{sec:gm}, the nested algebraic Bethe ansatz is used to tackle the diagonalisation of $X_{\tm}$.

\subsection{Examples}\label{subsec-exa}

\paragraph{1.} First assume that the multiplicity $d_{\tm}$ is equal to 1. This is not very interesting since there is no missing label in this case. Nevertheless, the results above 
provide the values of $X_{\tm}$ and $Y_{\tm}$. They are simply numbers (1 by 1 matrices) given, respectively, 
by the formulas for $a_{11}$ and $b_{11}$ of the preceding subsection. For example, assume that $d_{\tm}=\ell=1$. Then the value of $X_{\tm}$ is:
\begin{eqnarray}
 x&=&\frac{1}{6}\left(  l_{m_1,m_2} -l_{m'_1,m'_2} + l_{m''_1,m''_2}+(n-1)(2m'_1+4m'_2-3)(m_1+m'_1-n)-m_1+m_2-m'_1+m'_2 \right)\nonumber\\
 && -\frac{m'_1m'_2(3m_1-3m_2+m'_1-m'_2)}{27}  -\frac{(m_1+m_2)(2m'_1-1+2m'_2)(m'_1-m'-2)}{12}\nonumber\\
&&+\left( \frac{m'_1-m'_2}{27}+ \frac{m_1-m_2}{36} \right)(m'_1+m'_2)(2m'_1+2m'_2-9) \ .\label{xsca}
\end{eqnarray}

\paragraph{2.} We consider a situation where the multiplicity $d_{\tm}$ is equal to $2$, namely, we take:
\begin{equation}
 [m_1,m_2]=[2,2]\,,\ \ \ \ \ \ [m'_1,m'_2]=[m''_1,m''_2]=[p+1,q+1]\,, \label{ex2}
\end{equation}
where $p,q\geq 1$. In this case, we have $n=\ell=2$, and the representation $[m''_1,m''_2]$ appears with multiplicity 2 in the tensor product $[m_1,m_2]\otimes [m'_1,m'_2]$. In this example, we have:
\[X_{\tm}=\left(\begin{array}{cc}
\alpha-(p+q+\frac{3}{2}) & -p \\ q & \alpha+(p+q+\frac{3}{2})
\end{array}\right)\,,\ \ \ \ \ \text{where $\alpha=-\frac{1}{27}(p-q)(3+2p+q)(3+p+2q)$}\ .\]
 The two eigenvalues of $X_{\tm}$, which can serve as a missing label, are:
 \begin{equation}
  \alpha\pm\sqrt{\left(p+q+\frac{3}{2}\right)^2-p q}\, .   \label{exx2}
 \end{equation}
In this example, the matrix representing $Y$ is:
\[Y_{\tm}=-\frac{1}{12}\left(\begin{array}{cc}
9+24p+8p^2+4pq& 4p(6+p+2q) \\ 4q(6+2p+q) & 9+24p+8q^2+4pq
\end{array}\right)\ ,\]
and its eigenvalues, which can also serve as a missing label, are:
\[-\frac{3}{4}\ \ \ \text{and}\ \ \ \ -\frac{1}{12}(9+24p+24q+8p^2+8q^2+8pq)\ .\]

\paragraph{3.} Finally, let us consider the following situation:
\begin{equation}
 [m_1,m_2]=[2,p+1]\,,\ \ \ \ \ \ [m'_1,m'_2]=[2,q+1]\,,\ \ \ \ \ \ \ [m''_1,m''_2]=[L,p+q-2L+4]\,, \label{ex3}
\end{equation}
where $2\leq L\leq p+1\leq q+1$. In this case, we have $\ell=L$, $n=2$ and the multiplicity $d_{\tm}$ is equal to $2$. We can calculate explicitly the matrices representing $X_{\tm}$ and $Y_{\tm}$, and their eigenvalues. One finds for the eigenvalues of $X_{\tm}$:
\begin{eqnarray}
\frac{1}{54}(p-q)(54+36p+36q+4p^2+4q^2+10pq-36L+9L^2-9L p-9L q) \nonumber\\
\qquad \pm\sqrt{\frac{1}{4}(1-4L+L^2-L p-L q)^2-(2-L+q)(2-L+p)}\ .  \label{exx3}
\end{eqnarray}
Some formulas conjectured in \cite{CaSt} follow from this calculation. We note that for $L=2$, we recover the same eigenvalues than in the preceding examples. This is a particular example of the symmetries that we are going to study in Section \ref{sec-sym}.

\subsection{Clebsch--Gordan coefficients of $\mathfrak{su}_2$ and Heun--Hahn operators \label{sec:CGsu2}} 

In this subsection, a surprising relation between the operators appearing naturally in the computation of the Clebsch--Gordan coefficients of $\mathfrak{su}_2$  \cite{Alexei}
and the matrices $X_{\tm}$ and $Y_{\tm}$ given by \eqref{eq:XY} is pointed out.

Let us recall the definition of the Clebsch--Gordan coefficients of $\mathfrak{su}_2$.
The $\mathfrak{su}_2$ Lie algebra is generated by $s_x,s_y,s_z$ satisfying $[s_x,s_y]=is_z$, $[s_y,s_z]=is_x$ and $[s_z,s_x]=is_y$. 
Its finite irreducible representations are characterized by an half-integer $J\in \mathbb{Z}_{>0}/2$, called spin and the vectors of 
this representation are $|JM\rangle$ with $-J\leq M \leq J$ where $M$ is the eigenvalues of $s_z$.
There exist two natural bases when the tensor product of two finite representations of spin $J_1$ and $J_2$ is considered.
The vectors of the uncoupled basis are given by, for $-J_i\leq M_i \leq J_i$,
\begin{equation}
 |J_1 J_2;M_1 M_2\rangle = |J_1M_1\rangle \otimes |J_2M_2\rangle.
\end{equation}
The coupled vectors are $|J_1J_2;J_3M_3\rangle$ ( $|J_1-J_2|\leq J_3\leq J_1+J_2$ and $-J_3\leq M_3 \leq J_3$) 
where $M_3$ is the eigenvalue of $s^{(1)}_z +s^{(2)}_z$ ( $s^{(1)}_z=s_z\otimes 1$ and $s^{(2)}_z=1\otimes s_z$).
The expansion coefficients $\langle J_1 J_2; M_1 M_2 \ | \ J_1J_2;J_3M_3 \rangle$ are called 
the Clebsch--Gordan coefficients of $\mathfrak{su}_2$.

If $J_1$, $J_2$ and $M_3=M_1+M_2$ are fixed, the two remaining pertinent 
operators are $H_1=\frac{1}{2}(s^{(1)}_z-s^{(2)}_z)$ and $H_2=\sum_{a=x,y,z} (s^{(1)}_a+s^{(2)}_a)^2$.
In the uncoupled basis, $H_1$ is diagonal and $H_2$ is tridiagonal (see \cite{Alexei} for $H_1$ and $H_2$ in the coupled basis).
More precisely, 
the vectors of the uncoupled basis are given by
\[v_{i}=(-1)^i |J_1 J_2; (M_3+J_2-i) (i-J_2)\rangle\quad \text{with }\qquad \eta_b \leq i \leq \eta_a-1\]
where 
\[\eta_a=\text{min}\{\eta_1,\eta_2\}\quad \text{and}\qquad \eta_b=\text{max}\{\eta_5,\eta_6\}\ ,\]
and the parameters $\eta$ are 
\begin{subequations}\label{eq:eta}
\begin{eqnarray}
 &&(\eta_1\,,\ \eta_2\,,\ \eta_3)=(J_1+J_2+M_3+1\,,\ 2J_2+1\,,\ J_2+J_3+M_1+1   )\,,\\
 &&(\eta_4\,,\ \eta_5\,,\ \eta_6)=(J_2+J_3+1\,, J_2-J_1+M_3\,,\ 0   )\,.
\end{eqnarray}
\end{subequations}
In this uncoupled basis, $H_1$ becomes 
\[ H_1=\text{diag}\left(0,-1,\dots,-N+1  \right)+\left(\frac{\eta_1+\eta_2+\eta_5+\eta_6-2}{4}-\eta_b\right) \mathbf{1}_N, \]
with $N=\eta_a-\eta_b$ and $H_2$ is the following tridiagonal matrix
\begin{equation}\label{H2}
 H_2=\left(\begin{array}{ccccc}
    \alpha_{11} &  \alpha_{12}  &0  & \dots  & 0  \\
     \alpha_{21} &  \alpha_{22}  & \ddots  & \ddots & \vdots  \\
     0 & \ddots  & \ddots & \ddots & 0 \\
     \vdots & \ddots & \ddots & \ddots &  \alpha_{N-1,N}\\
     0 & \dots & 0 &     \alpha_{N,N-1}  &   \alpha_{N,N}
 \end{array}\right)
\end{equation}
where the coefficients are given by:
\[\begin{array}{r}\alpha_{j,j+1}=-(j+\eta_b-\eta_1)(j+\eta_b-\eta_2)\ ,\qquad 
\alpha_{j+1,j}=-(j+\eta_b-\eta_5)(j+\eta_b-\eta_6)\ ,
\end{array}\]
\[\alpha_{j,j}=\alpha_{j,j-1}+\alpha_{j,j+1} +\frac{1}{4}(\eta_1+\eta_2-\eta_5-\eta_6)(\eta_1+\eta_2-\eta_5-\eta_6-2) \ .\]
The matrices $H_1$ and $H_2$ are called Hahn matrices since they appear in the study of the Hahn polynomials. 
Indeed their entries are the recurrence coefficients of the Hahn polynomials \cite{KLS}.

The connection of these matrices with $X_{\tm}$ and $Y_{\tm}$ given by \eqref{eq:XY} is provided by the following proposition.
\begin{prop} \label{lem1}
Let us replace the elements of $\{\eta_1\,,\ \eta_2\,,\ \eta_3\}$ (resp.  $\{\eta_4\,,\ \eta_5\}$)  by the elements of $\{\xi_1\,,\ \xi_2\,,\ \xi_3\}$ (resp.  $\{\xi_4\,,\ \xi_5\}$) 
such that  $\eta_a$ is replaced by $\xi_a$ (resp. $\eta_b$ is replaced by $\xi_b$) in the matrices $H_1$ and $H_2$. We recall that $\xi_6=\eta_6=0$.
In particular, one gets $N=d_\tm$.
The matrices $X_{\tm}$ and $Y_{\tm}$ are expressed in terms of the Hahn matrices $H_1$ 
and $H_2$ as follows
\begin{subequations}
 \begin{eqnarray}
 X_{\tm}&=& z_0+ z_1 H_1 +z_2 H_2 +z_3 H_3 +z_4 \{H_1,H_2\}\ , \label{eq:Xe} \\
 Y_{\tm}&=& z_5+z_6 H_1+z_7 H_2+z_8 H_3+z_9 \{H_1,H_2\} +z_{10} H_1^2+z_{11} H_1^2 H_2+z_{12} H_1 H_2H_1 \ , \label{eq:Ye}
\end{eqnarray}
\end{subequations}
with the parameters $z$ given in Appendix \ref{app:B}.
\end{prop}
\proof  From expressions \eqref{eq:Xe} of $X_{\tm}$, one computes its non-vanishing entries
\begin{subequations}
\begin{eqnarray}
 (X_{\tm})_{i,i+1}&=&(z_2+z_3-(2i-1)z_4) \alpha_{i,i+1}\ , \\
 (X_{\tm})_{i+1,i}&=&(z_2-z_3-(2i-1)z_4) \alpha_{i+1,i}\ , \label{eq:cXd}\\
 (X_{\tm})_{i,i}&=&z_0+z_2\alpha_{i,i}   -\left(i-\frac{\eta_1+\eta_2+\eta_4+\eta_5+2}{4}+\eta_b\right)(z_1+2z_4 \alpha_{i,i})\,.
\end{eqnarray}
\end{subequations}
By using the explicit forms of the parameters $z$ given in Appendix \ref{app:B}, we recover the expressions \eqref{eq:XY} for $X_\tm$.
The proof is similar for $Y_\tm$.
 \endproof

\begin{remax}
 The quadratic combination of Hahn matrices given in the R.H.S. of relation \eqref{eq:Xe} is usually called a algebraic Heun--Hahn operator.
 The notion of algebraic Heun operator which generalizes the standard Heun operator were recently introduced in \cite{GVZ} and has been studied
 in different cases \cite{VZ,BTVZ,CVZ,BCTVZ}.
 \end{remax}

\subsection{Bethe ansatz to diagonalize $X_\tm$ \label{sec:gm}}

As mentioned previously, there does not exist, in general, analytical formulas for the eigenvalues of $X_\tm$ and $Y_\tm$.
However, the connection with a Heun--Hahn operator suggests that it is possible to diagonalize $X_\tm$ with the help of the Bethe ansatz 
as it has been done in \cite{BP, BCSV} for different algebraic Heun operators. The Heun--Hahn operator $X_\tm$ has been identified
in the context of the two-body homogeneous rational Gaudin models in \cite{CAM} and the nested Bethe ansatz developed in \cite{Ju} for these models 
allows to obtain the eigenvalues of $X_\tm$ in terms of Bethe roots. Let us remark that the application of a similar approach to the diagonalisation of $Y_\tm$ is an open problem since
this operator does not appear in the context of integrable systems.

The details of the implementations of the nested Bethe ansatz associated to the diagonalization of $X_\tm$ are not 
given here but the results displayed below have been obtained by following the lines of \cite{Ju}.
In this context the eigenvalues of $X_\tm$ are associated to two sets of complex numbers
$\bar \nu=\{\nu_1,\dots \nu_{n-1}\}$ and $\bar \lambda=\{\lambda_1,\dots \lambda_{\ell-1}\}$ (where $\ell$ and $n$ are given by \eqref{eq:nl2}) 
satisfying the following equations 
\begin{eqnarray}
 &&\frac{m_1-1+\nu_p(2+m'_1-2n+\ell)}{\nu_p(1-\nu_p)}
 = \sum_{\genfrac{}{}{0pt}{2}{k=1}{k\neq p}}^{n-1} \frac{2}{\nu_p-\nu_k} -\sum_{r=1}^{\ell-1} \frac{1}{\nu_p-\lambda_r} \quad \text{for} \quad 1\leq p \leq n-1 \ , \label{eq:BE31q}  \\
 &&\frac{m_2-1+\lambda_s(2+m'_2-2\ell+n)}{\lambda_s(1-\lambda_s)}= \sum_{\genfrac{}{}{0pt}{2}{r=1}{r\neq s}}^{\ell-1} \frac{2}{\lambda_s-\lambda_r} 
 -\sum_{k=1}^{n-1} \frac{1}{\lambda_s-\nu_k} \quad \text{for} \quad 1\leq s \leq \ell-1 \ . \label{eq:BE32q} 
\end{eqnarray}
These equations are called Bethe equations and their solutions, Bethe roots. The physical Bethe roots correspond to 
the Bethe roots pairwise distinct\footnote{There exists exceptional points where two Bethe roots must be equal to obtain 
the complete spectrum. We will not discuss these cases here.} and all different from $0$ and $1$.
For each physical Bethe roots, the eigenvalue of $X_\tm$ is given by
\begin{eqnarray}
&& x-\frac{\ell-1}{6} \left(  (m_1-m_2+2m'_1-2m'_2 -3n+3\ell )(n-1)+(m_2+2m_1-2m'_1-m'_2)(\ell-m_2-m'_2 ) \right) \nonumber
\\
&&+\frac{m_1+m_2-1}{2}\sum_{k=1}^{n-1}\sum_{s=1}^{\ell-1}\frac{\nu_k+\lambda_s-2}{\nu_k-\lambda_s} \, ,\label{eq:x}
\end{eqnarray}
where $x$ is given by \eqref{xsca}.
Let us also mention that the nested Bethe ansatz provides the eigenvectors even though we do not discuss this point here.

Instead of solving directly the Bethe equations, it is more convenient to solve differential equations 
whose zeros of the polynomial solutions provide the Bethe roots.
Following \cite{MV,MV2,MV3}, one gets the following result.
The parameters $\bar \nu=\{\nu_1,\dots \nu_{n-1}\}$ and $\bar \lambda=\{\lambda_1,\dots \lambda_{\ell-1}\}$, pairwise distinct, 
are solutions of the Bethe equations \eqref{eq:BE31q} and \eqref{eq:BE32q}
if and only if there exist two polynomials $\rho_1(u)$ and $\rho_2(u)$ such that the coupled differential equations
\begin{eqnarray}
 u(u-1)y_2(u)y_1''(u)+[ ( u(m_1'-2n+\ell+2)+m_1-1)y_2(u)+u(1-u)y_2'(u) ]y_1'(u)+\rho_1(u)y_1(u)=0\,,\label{eq:PJ1}\\
 u(u-1)y_1(u)y_2''(u)+[ ( u(m_2'-2\ell+n+2)+m_2-1)y_1(u)+u(1-u)y_1'(u) ]y_2'(u)+\rho_2(u)y_2(u)=0\,, \label{eq:PJ2}
\end{eqnarray}
have polynomial solutions $y_1(u)$ and $y_2(u)$ of degree $n-1$ and $\ell-1$, respectively. 
Then, the Bethe roots are the zeros of these polynomials: $\displaystyle y_1(u)=\prod_{k=1}^{n-1} (u-\nu_k)$ and $\displaystyle y_2(u)=\prod_{r=1}^{\ell-1} (u-\lambda_r)$.
The polynomial $\rho_1(u)$ (resp. $\rho_2(u)$) has degree $\ell-1$ (resp. $n-1$) and the coefficient in front of $u^{\ell-1}$ (resp. $u^{n-1}$) is $(n-1)(n-1-m'_1)$ (resp. $(\ell-1)(\ell-1-m'_2)$), computed by 
taking the asymptotic $u\to\infty$ in \eqref{eq:PJ1} (resp. \eqref{eq:PJ2}). These equations are closely related to the Jacobi-Pineiro polynomials \cite{MV2}.

In \cite{MV2}, a third equivalent way to provide the Bethe equations is also given. In the case treated here, it says that 
the parameters $\{\nu_p\ |\ p=1, \dots n\}$ and  $\{\lambda_r\ |\ r=1, \dots \ell\}$, pairwise distinct, are solutions of the Bethe equations \eqref{eq:BE31q}-\eqref{eq:BE32q}
if and only if there exist two polynomials $\displaystyle y_1(u)=\prod_{k=1}^n (u-\nu_k)$ and $\displaystyle y_2(u)=\prod_{r=1}^\ell (u-\lambda_r)$ with simple roots, solutions of
\begin{eqnarray}
 && u(u-1)\big( y_1(u)y_2''(u)  - y_1'(u)y_2'(u)+y_1''(u)y_2(u) \big) \nonumber \\
 &&+ \big( u(m_1'-2n+\ell+2)+m_1-1\big)y'_1(u)y_2(u) + \big(u(m_2'-2\ell+n+2)+m_2-1\big)y_1(u)y'_2(u)\nonumber\\
 &&+\big( (\ell-1)(\ell-1-m_2')+(n-1)(n-1-m_1') -(n-1)(\ell-1)\big)y_1(u)y_2(u)=0\, . \label{eq:hh}
\end{eqnarray}

The problem to prove that there are exactly $d_\tm$ different polynomial solutions to the Bethe equations \eqref{eq:BE31q}-\eqref{eq:BE32q} (or \eqref{eq:PJ1}-\eqref{eq:PJ2})
is a difficult task (and it is even wrong in a few cases see \cite{MV2}).
In the following, the solution to this problem is given for the examples presented in Section \ref{subsec-exa}. The Bethe ansatz reproduces the results obtained by direct diagonalization of $X_\tm$.

\paragraph{Example 1.} In this example, $\ell=1$ and $d_\tm=1$. In this case, the differential equations \eqref{eq:PJ1}-\eqref{eq:PJ2} reduce to ($y_2(u)=1$)
\begin{eqnarray}
u(u-1)y_1''(u)+( u(m_1'-2n+3)+m_1-1)y_1'(u)+(n-1)(n-1-m_1') y_1(u)=0\ .\label{eq:BE2}
\end{eqnarray}
The Jacobi polynomial \cite{KLS} (the normalisation is fixed by requiring that the polynomial be monic)
 \begin{equation}
 y_1(u)=\frac{(n-1)!(m_1'-n)!}{(m'_1-1)!}P_n^{(m_1+m_1'-2n+1,-m_1)}(2u-1) \ , \label{eq:jac}
 \end{equation}
is a polynomial solution to \eqref{eq:BE2}.
The Bethe ansatz provides exactly one solution and, setting $\ell=1$ in \eqref{eq:x}, we recover the unique eigenvalue of $X_\tm$ given by \eqref{xsca}.

\paragraph{Example 2.} Let us consider the case given by \eqref{ex2}. Here $n=\ell=2$ and the Bethe equations \eqref{eq:BE31q}-\eqref{eq:BE32q} become
\begin{eqnarray}
1+\nu p = -\frac{\nu(1-\lambda)}{\nu-\lambda} \quad \text{and} \qquad 1+\lambda q=-\frac{\lambda(1-\nu)}{\lambda-\nu} \ .  \label{eq:BEnl_{1}}
\end{eqnarray}
To simplify the notation, we replace $\nu_1$ by $\nu$ and $\lambda_1$ by $\lambda$.
There exist two solutions to these equations given by
\begin{eqnarray}
 \lambda = \frac{-2p-4q-3\pm 2\sqrt{\left(p+q+\frac{3}{2}\right)^2-p q}}{2q(p+q+1)}\quad , \qquad
 \nu= \frac{-4p-2q-3\pm 2\sqrt{\left(p+q+\frac{3}{2}\right)^2-p q}}{2p(p+q+1)}\, .
\end{eqnarray}
Putting these values of Bethe roots in \eqref{eq:x}, one recovers the explicit expression \eqref{exx2} given previously.

\paragraph{Example 3.} Finally, let us consider the case given by \eqref{ex3}. 
The Bethe equations provided by \eqref{eq:hh} reduce to ($y_1(u)=u-\nu$)
\begin{eqnarray}
 &&u(u-1)(u-\nu)y_2''(u)+[ ( u(q-2\ell+5)+p)(u-\nu)+u(1-u) ]y_2'(u)\nonumber\\
 && \qquad+((\ell-1)(\ell-2-q)(u-\nu)+\nu \ell+1)y_2(u)=0.  \label{eq:y2s}
\end{eqnarray}
Again, we have used $\nu$ instead of $\nu_1$. 
The following combinations of Jacobi polynomials 
\begin{equation}
 y_2(u)=\frac{(\ell-1)!(q-\ell+1)!}{q!  }\left[P_{\ell-1}^{(q+p-2\ell+4,-p-2)}(2u-1)  +\frac{\nu+1}{\ell-1}P_{\ell-2}^{(q+p-2\ell+4,-p-1)}(2u-1)  \right]\ , \label{eq:yy2}
\end{equation}
satisfies relation \eqref{eq:y2s} if $\nu$ is one of the zeros of the following quadratic polynomial
\begin{equation}
 (q+2)\nu^2-( \ell^2-\ell(p+q+4)-1)\nu+p+2=0 \ . \label{eq:nuz}
\end{equation}
This result is proven using the second order differential equation satisfied by the Jacobi polynomials 
\eqref{eq:BE2} and the properties
\begin{eqnarray}
 (u-1) \frac{d}{du}P_\ell^{(q+p-2\ell+2,-p-2)}(2u-1)&=&\ell P_\ell^{(q+p-2\ell+2,-p-2)}(2u-1)\nonumber\\
&&-(q+p-\ell+2)P_{\ell-1}^{(q+p-2\ell+2,-p-1)}(2u-1)\ ,\\
 u(u-1)\frac{d}{du}P_{\ell-1}^{(q+p-2\ell+2,-p-1)}(2u-1)&=& \ell  P_\ell^{(q+p-2\ell+2,-p-2)}(2u-1)\nonumber \\
 &&- ( u(q+1-\ell)+p+1) P_{\ell-1}^{(q+p-2\ell+2,-p-1)}(2u-1) \ .
\end{eqnarray}
The two corresponding eigenvalues of the operators $X$ are given by \eqref{exx3}.

\section{Symmetry}\label{sec-sym}

Given two sets of parameters $\tm$ and $\tm'$, it is of interest to determine when the two pairs of matrices $(X_{\tm},Y_{\tm})$ and $(X_{\tm'},Y_{\tm'})$ are conjugated 
(we mean simultaneously conjugated by the same change of basis). If this is so the eigenvalues of $X$ and of $Y$ (the possible missing labels) will be invariant when $\tm$ is replaced by $\tm'$.
Since $X$ and $Y$ generate, up to the central elements, the whole centraliser, any missing label operator taken from the centraliser $Z_2(\mathfrak{su}_3)$ will have the same eigenvalues in the physical representation associated to $\tm$ and in the physical representation associated to $\tm'$. 

If $(X_{\tm},Y_{\tm})$ and $(X_{\tm'},Y_{\tm'})$ are conjugated, we call the corresponding transformation of the parameters $\tm$ into $\tm'$ a ``symmetry'' of the missing label operators (in fact we will consider the slightly more general situation where $X_{\tm}$ and $X_{\tm'}$ are conjugated up to a sign).

In this section, we exhibit a group of $144$ symmetries of the missing label operators. The natural symmetries (permuting the representations and taking the dual) 
form a group of order 12, so we have many interesting additional symmetries. 
It turns out that these $144$ symmetries are exactly the symmetries of the Littlewood--Richardson $d_{\tm}$ described in \cite{BR}.
We note the remarkable fact that all the symmetries of  $d_{\tm}$ actually extend to symmetries of the operators $X_{\tm}$ and $Y_{\tm}$.

We give below two graphical representations of the symmetry group of order 144. One of these interpretations is with magic squares of numbers, 
which is similar to the ones of the Clebsch--Gordan coefficients of $\mathfrak{su}_2$. 
The other one uses affine transformations on a torus.
We will provide also an interpretation of the group of $144$ symmetries in terms of the root system of type $E_6$ in Section \ref{sec-alg}.

\subsection{Symmetry of the missing label operators $X_{\tm}$ and $Y_{\tm}$}

\paragraph{Magic squares.} Say $n\leq \ell$. Then the arrangement of the Theorem stated in the introduction is given in terms of the parameters $\xi_1,\dots,\xi_5$ (from the beginning of Section \ref{subsec-formulas}) by:
\begin{equation}\label{square-xi}\left[\begin{array}{ccc|ccc}
\xi_1-\xi_5 & \xi_2 & \xi_3-\xi_4                             \ \ \ & \ \ \ \\[0.5em]
\xi_2-\xi_4 & \xi_3-\xi_5 & \xi_1      \ \ \  & \ \ \ & +(\ell-n) &  \\[0.5em]
\xi_3 & \xi_1-\xi_4 & \xi_2-\xi_5 \ \ \ & \ \ \
\end{array}
\right]
\end{equation}
where the right square is equal to the left square with $(\ell-n)$ added to all entries. If $\ell\leq n$ then the arrangement of the Theorem  in the introduction is the same with the two squares exchanged. So up to exchanging the left and right squares, we can assume that the arrangement is as in (\ref{square-xi}).

Note that the parameters $\xi_1,\dots,\xi_5$ and $\ell-n$ uniquely determine the parameters in $\tm$, so that the transformations described in the theorem below uniquely determine a transformation on the parameters $\tm$.
\begin{theo}\label{theo-sym}
The pair of matrices $(X_{\tm},Y_{\tm})$ is conjugated to the pair of matrices $(\pm X_{\tm'},Y_{\tm'})$ for a transformation of the parameters $\tm\mapsto \tm'$ corresponding to the following operations:
\begin{itemize}
\item permutations of the lines and columns simultaneously on the left and right $3\times 3$ squares;
\item transposition simultaneously on the left and right $3\times 3$ squares;
\item exchange of the left and right $3\times 3$ squares.
\end{itemize}
The sign in front of $X_{\tm'}$ is $(-1)$ to the power the number of transpositions of lines and columns.
\end{theo}
The proof of this Theorem is postponed to Section \ref{sec:pr}.

Let us explain how these transformations generate a group of order 144, and what is the structure of this group. 
For a permutation $\pi$ in the symmetric group $S_3$, $\pi_{L}$ (resp. $\pi_{C}$) denotes the corresponding permutation of lines (resp. columns) of the two $3\times 3$ squares. Let also $T$ denote the simultaneous transposition and $S$ the exchange of the two squares.
The symmetry $S$ commutes with all the other transformations, any $\pi_C$ commutes with any $\pi'_L$ and $T\circ \pi_L=\pi_C\circ T$ for all $\pi\in S_3$. 
It follows that any element of the group generated by the transformations in the theorem above can be written in the following form:
\[\pi_C\circ\pi'_L\circ T^{\epsilon_1}\circ S^{\epsilon_2}\,,\ \ \ \ \ \text{where $\pi,\pi'\in S_3$ and $\epsilon_1,\epsilon_2\in\{0,1\}$.}\]
This explains the order $144=6\times 6\times 2\times 2$ of the group. The structure of the group is a semi-direct product as follows:
\begin{equation}\label{group1}\Bigl((S_3\times S_3)\rtimes \mathbb{Z}/2\mathbb{Z}\Bigr)\times \mathbb{Z}/2\mathbb{Z}\ .
\end{equation}

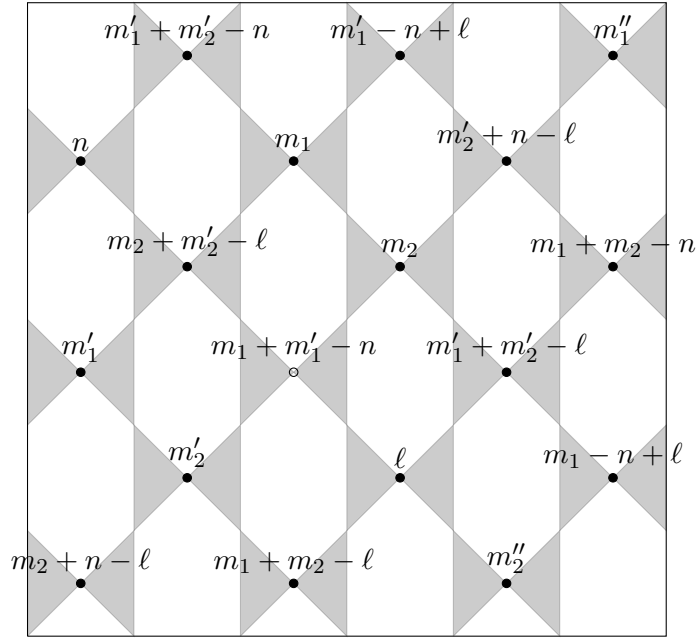
\begin{figure}[htb]
\begin{center}
\begin{tikzpicture}[scale=0.7]
\draw[fill] (1,1) circle [radius=0.08];
\node [above] at (1,1) {$m_2+n-\ell$};
\draw[fill] (1,5) circle [radius=0.08];
\node [above] at (1,5) {$m'_1$};
\draw[fill] (1,9) circle [radius=0.08];
\node [above] at (1,9) {$n$};
\draw[fill] (3,11) circle [radius=0.08];
\node [above] at (3,11) {$m'_1+m'_2-n$};
\draw[fill] (3,7) circle [radius=0.08];
\node [above] at (3,7) {$m_2+m'_2-\ell$};
\draw[fill] (3,3) circle [radius=0.08];
\node [above] at (3,3) {$m'_2$};

\draw[fill] (5,1) circle [radius=0.08];
\node [above] at (5,1) {$m_1+m_2-\ell$};
\draw (5,5) circle [radius=0.08];
\node [above] at (5,5) {$m_1+m'_1-n$};
\draw[fill] (5,9) circle [radius=0.08];
\node [above] at (5,9) {$m_1$};
\draw[fill] (7,11) circle [radius=0.08];
\node [above] at (7,11) {$m'_1-n+\ell$};
\draw[fill] (7,7) circle [radius=0.08];
\node [above] at (7,7) {$m_2$};
\draw[fill] (7,3) circle [radius=0.08];
\node [above] at (7,3) {$\ell$};

\draw[fill] (9,1) circle [radius=0.08];
\node [above] at (9,1) {$m''_2$};
\draw[fill] (9,5) circle [radius=0.08];
\node [above] at (9,5) {$m'_1+m'_2-\ell$};
\draw[fill] (9,9) circle [radius=0.08];
\node [above] at (9,9) {$m'_2+n-\ell$};
\draw[fill] (11,11) circle [radius=0.08];
\node [above] at (11,11) {$m''_1$};
\draw[fill] (11,7) circle [radius=0.08];
\node [above] at (11,7) {$m_1+m_2-n$};
\draw[fill] (11,3) circle [radius=0.08];
\node [above] at (11,3) {$m_1-n+\ell$};
\draw[thin] (0,0)--(0,12)--(12,12)--(12,0)--(0,0);

\foreach \y in {0,4,8}
\foreach \x in {0,4,8} 
{\draw [fill, opacity=0.2] (\x,\y)--(\x+1,\y+1)--(\x,\y+2)--(\x,\y);
\draw [fill, opacity=0.2] (\x+1,\y+1)--(\x+2,\y+2)--(\x+2,\y)--(\x+1,\y+1);
\draw [fill, opacity=0.2] (\x+2,\y+2)--(\x+3,\y+3)--(\x+2,\y+4)--(\x+2,\y+2);
\draw [fill, opacity=0.2] (\x+3,\y+3)--(\x+4,\y+4)--(\x+4,\y+2)--(\x+3,\y+3);
}

\end{tikzpicture}

\caption{Symmetries of the missing labels seen as affine transformations on a torus. \label{fig}}
\end{center}
\end{figure}

\paragraph{Affine transformations on a torus.} Another way to describe this symmetry group consists in drawing the $18$ elements of the set \eqref{eq:deg} on a torus as displayed in Figure \ref{fig}.  
The opposite sides of the square are identified.
The $144$ transformations of Theorem \ref{theo-sym} correspond to the affine transformation leaving stable the set of dots of this picture. 
Given one of the 18 dots, the transformations fixing this dot are the symmetries of a square of which this dot is the center (a dihedral group of order 8). 
The full symmetry is the semi-direct product of one of these dihedral groups with the subgroup of translations leaving the set of dots invariants. These translations are generated, for example, by two translations of order 6: say, $t_1$ sending $m_1$ to $m_2$, and $t_2$ sending $m_1$ to $m_2+m'_2-\ell$. These two translations satisfy $t_1^3=t_2^3$, so that the translation subgroup is of order $18$. This checks with $144=8\times 18$.

Moreover, if a symmetry maps gray zones to gray zones then it corresponds to a transformation $\tm\mapsto \tm'$ such that $X_{\tm}$ and $X_{\tm'}$ are conjugated. Otherwise, it corresponds to a case where $X_{\tm}$ and $-X_{\tm'}$ are conjugated.

Let us emphasize that these two graphical representations for the symmetries of the missing label $X_\tm$ and $Y_\tm$ provide also new interpretations for the symmetries of the Littlewood--Richardson coefficients 
found in \cite{BR}.

\paragraph{Classical symmetries.} There are 12 classical symmetries generated by the permutations of the three representations involved and the replacement of all three representations by their dual. 

The permutations of the three representations correspond naturally to the 6 permutations of the set $\{(m_1,m_2),\,(m'_1,m'_2),\ (m''_2,m''_1)\}$. 
Note the exchange of $m''_1,m''_2$: we need to take the dual of the third representation such that the three representations to be on an equal footing. 
These six transformations correspond to the (simultaneous) permutations of the columns of the squares in (\ref{square-xi}). 

Replacing all representations by their duals amounts to the following transformation on the parameters $m_1\leftrightarrow m_2$, $m'_1\leftrightarrow m'_2$ and $m''_1\leftrightarrow m''_2$. It corresponds to the exchange of the left and right $3\times 3$ squares composed with the exchange of the first and third lines.

\begin{remax}
The theorem has a very practical application. Given a set of parameters, 
the symmetry can be used to obtain a new set of parameters such that $\ell$ is the minimum of the 18 numbers in the arrangement. 
In this case, the formulas given in Section \ref{sec-XY} are easier to apply, since the parameters $\xi$'s are always the same, and the parameter $\xi_b$ is $0$.
\end{remax}

\paragraph{Examples.} Take the integers $L,p,q$ such that $2\leq L\leq p+1\leq q+1$. We consider the following situation:
\begin{equation}\label{exa1}
[m_1,m_2]=[2,p+1]\,,\ \ \ \ \ \ [m'_1,m'_2]=[2,q+1]\,,\ \ \ \ \ \ \ [m''_1,m''_2]=[L,p+q-2L+4]\,.
\end{equation}
The eigenvalues of $X_{\tm}$ are calculated in Example 3 in Section \ref{subsec-exa}. 
The corresponding arrangement of the two magic squares of the Introduction or of equation (\ref{square-xi}), is:
\[\left[\begin{array}{ccc|ccc}
2 & 2 & p+q-2L+4                             \ \ \ & \ \ \ L & L & p+q-L+2  \\[0.5em]
q+3-L & p+3-L & 2      \ \ \  & \ \ \ q+1 & p+1 & L \\[0.5em]
p+3-L & q+3-L & 2 \ \ \ & \ \ \ p+1 & q+1 & L
\end{array}\right]\]
Now exchange the first two lines (in the two squares). This gives:
\[\left[\begin{array}{ccc|ccc}
q+3-L & p+3-L & 2                             \ \ \ & \ \ \ q+1 & p+1 & L  \\[0.5em]
2 & 2 & p+q-2L+4       \ \ \  & \ \ \  L & L & p+q-L+2 \\[0.5em]
p+3-L & q+3-L & 2 \ \ \ & \ \ \ p+1 & q+1 & L
\end{array}\right]\,,\]
which corresponds to the situation
\[[m_1,m_2]=[q+3-L,p+1]\,,\ \ \ \ \ \ [m'_1,m'_2]=[p+3-L,q+1]\,,\ \ \ \ \ \ \ [m''_1,m''_2]=[L,2]\,.\]
Here we have $\ell=p+q-L+2$ and $n=p+q-2L+4$, and we find that $n\leq \ell$ and $d_{\tm}=2$. 
If one uses the explicit formulas \eqref{eq:XY} for $X_{\tm}$ noting that $\xi_b=p+1-L$, its eigenvalues are the opposite of the eigenvalues for the previous case. 
The theorem of this section indicates that this must be so, without any calculation.

Finally, in the last arrangement, perform the exchange of the first and third columns simultaneously in the two squares. This gives
\[\left[\begin{array}{ccc|ccc}
2 & p+3-L & q+3-L                           \ \ \ & \ \ \ L & p+1 & q+1  \\[0.5em]
 p+q-2L+4  & 2 & 2      \ \ \  & \ \ \  p+q-L+2 & L & L \\[0.5em]
2 & q+3-L & p+3-L \ \ \ & \ \ \ L & q+1 & p+1
\end{array}\right]\,,\]
which corresponds to the situation
\[[m_1,m_2]=[2,L]\,,\ \ \ \ \ \ [m'_1,m'_2]=[p+3-L,q+1]\,,\ \ \ \ \ \ \ [m''_1,m''_2]=[p+1,q+3-L]\,.\]
Again the direct computation of the eigenvalues of $X_{\tm}$ shows that they are the same as for the parameters (\ref{exa1}), as ensured by the theorem. For $L=2$, 
this explains the coincidence observed between the examples 2 and 3 in Section \ref{subsec-exa}.

\subsection{Proof of the theorem \label{sec:pr}} 

The proof of the Theorem is given for the different types of transformations.

\paragraph{Exchange of the two squares.} It is simpler in this case to look directly at the squares given in the introduction. The transformation does the following:
\[n\leftrightarrow\ell\,,\ \quad m_2\leftrightarrow m_2+n-\ell\,,\ \quad m'_1+\ell-n\leftrightarrow m'_1\ \quad\text{and}\ \quad \ell-m'_2\mapsto \ell-m'_2\,,\ \quad n-m_1\mapsto n-m_1\,.\]
Since it exchanges $n$ and $\ell$, it also exchanges one of the inequalities ($n\leq \ell$ or $\ell\leq n$) with the other. 
Looking at the definition of the parameters $\xi$'s in each case, we conclude that they are transformed according to:
\[\{\xi_1,\xi_2,\xi_3\}\mapsto\{\xi_1,\xi_2,\xi_3\}\ \ \quad\ \ \{\xi_4,\xi_5,\xi_6\}\mapsto\{\xi_4,\xi_5,\xi_6\}\ .\]  So the size of the matrices stays the same (\emph{i.e.} 
this is a symmetry of the Littlewood-Richardson coefficient). Then we note that the parameters $\Lambda$ and $\lambda_{\pm}$ are fixed. So the matrices $X_{\tm'}$ and $Y_{\tm'}$ 
are equal to, respectively, $X_{\tm}$ and $Y_{\tm}$.\\

\paragraph{Simultaneous transposition into the two squares.} 
Let us remark that the simultaneous operations on the left and right squares fix the difference $n-\ell$, and therefore are completely given once their actions on the parameters $\xi$'s is given. 
The transposition into the squares is given on the parameters $\xi$'s by:
\[\left\{\begin{array}{l} \xi_1\mapsto \xi_1-\xi_4\ ,\\
\xi_2\mapsto \xi_2-\xi_4\ ,\\
\xi_3\mapsto \xi_3-\xi_4\ ,
\end{array}\right.\ \ \ \ \ \ \ \ \text{and}\ \ \ \ \ \ \ \ \left\{\begin{array}{l} \xi_4\mapsto -\xi_4\ ,\\
\xi_5\mapsto \xi_5-\xi_4\ .
\end{array}\right.\]
First, the minimum $\xi_a$ of $\xi_1,\xi_2,\xi_3$ and the maximum $\xi_b$ of $\xi_4,\xi_5,\xi_6$ are transformed as follows:
\[\xi_a\mapsto \xi_a-\xi_4\ \ \ \ \text{and}\ \ \ \ \ \xi_b\mapsto\xi_b-\xi_4\ .\]
The size of the matrices stays the same (\emph{i.e.} this is a symmetry of the Littlewood-Richardson coefficients). 
Then we note that the parameter $\Lambda$ is fixed while $\lambda_{\pm}$ is replaced by $\lambda_{\pm}-\xi_4$. So it follows that the off-diagonal coefficients 
of the new matrices coincide with those of the old ones. Next, we note that the transformation preserves the set $\{\xi_i-j-\xi_b+\frac{1}{2}\}_{i=1,\dots,6}$. 
Thus, since the diagonal coefficients of $X_{\tm}$ and $Y_{\tm}$ were expressed in terms of symmetric polynomials in this set, we conclude at once that the diagonal 
coefficients stay the same as well. So we conclude that the matrices $X_{\tm'}$ and $Y_{\tm'}$ are equal to, respectively, $X_{\tm}$ and $Y_{\tm}$.

\paragraph{Simultaneous permutations of columns or lines.}  Let us consider the permutation of the first two columns which corresponds to the following transformation:
\begin{equation}\label{transf-C12}
\left\{\begin{array}{l} \xi_1\mapsto \xi_1\,,\\
\xi_2\mapsto \xi_1-\xi_5\,,\\
\xi_3\mapsto \xi_1-\xi_4\,,
\end{array}\right.\ \ \ \ \ \ \ \ \text{and}\ \ \ \ \ \ \ \ \left\{\begin{array}{l} \xi_4\mapsto \xi_1-\xi_3\,,\\
\xi_5\mapsto \xi_1-\xi_2\,.
\end{array}\right.\,
\end{equation}
It is easy to see that all permutations of lines and of columns of the squares in (\ref{square-xi}) are of this form, 
up to permutations of $\{\xi_1,\xi_2,\xi_3\}$ and permutations of $\{\xi_4,\xi_5\}$. 
The verification is exactly the same for all these transformations (in other words, permutations of $\{\xi_1,\xi_2,\xi_3\}$ and of $\{\xi_4,\xi_5\}$ are obvious symmetries). 
Therefore it is enough to consider (\ref{transf-C12}).

Under the above transformation, the minimum $\xi_a$ of $\xi_1,\xi_2,\xi_3$ and the maximum $\xi_b$ of $\xi_4,\xi_5,\xi_6$ are transformed as follows:
\[\xi_a\mapsto \xi_1-\xi_b\ \ \ \ \text{and}\ \ \ \ \ \xi_b\mapsto\xi_1-\xi_a\ .\]
So the size of the matrices stays the same (\emph{i.e.} this is a symmetry of the Littlewood-Richardson coefficient). 
The parameter $\Lambda$ is fixed while $\lambda_{\pm}$ is replaced by $\xi_1-\lambda_{\mp}$ by this transformation. So the off-diagonal coefficients of the new matrices are:
\[\begin{array}{rcl}a'_{j,j+1}=(j-\xi_a+\xi_4)(j-\xi_a+\xi_5)(j-\xi_a+\xi_6) \quad & \text{and} & \quad b'_{j,j+1}=a'_{j,j+1}(j-\xi_a+\lambda_+)\ ,\\[0.5em]
a'_{j+1,j}=(j-\xi_a+\xi_1)(j-\xi_a+\xi_2)(j-\xi_a+\xi_3) \quad & \text{and} & \quad b'_{j+1,j}=a'_{j+1,j}(j-\xi_a+\lambda_-)\ .
\end{array}\]
We perform a complete inversion of the basis vectors, that is, we rename the indices setting $j'=\xi_a-\xi_b+1-j$. This allows to recover the off-diagonal coefficients of $X_{\tm}$ and $Y_{\tm}$ (with a minus sign for $X_{\tm}$). 
As for the diagonal coefficients, the transformation on the parameters together with the renaming of the indices transform the set $\{\xi_i-j-\xi_b+\frac{1}{2}\}_{i=1,\dots,6}$ to minus itself. Since the diagonal coefficients were expressed in terms of symmetric polynomials in this set, we recover the diagonal coefficients of $X_{\tm}$ and $Y_{\tm}$ (with a minus sign for $X_{\tm}$).

In conclusion, for any transformation corresponding to exchanging two lines or two columns simultaneously in the two squares in (\ref{square-xi}), the matrices $X_{\tm'}$ and $Y_{\tm'}$ are equivalent to, respectively, $-X_{\tm}$ and $Y_{\tm}$.

\begin{remax}
In fact, the $72$ transformations which act simultaneously on the two squares in (\ref{square-xi})  can be given explicitly. Choose $x,y,z$ such that $\{x,y,z\}=\{1,2,3\}$, and $s,t,u$ such that $\{s,t,u\}=\{4,5,6\}$. There are $36$ possible choices. Then set:
\[\left\{\begin{array}{l} \xi_1\mapsto \xi_x-\xi_s\\
\xi_2\mapsto \xi_y-\xi_{s}\\
\xi_3\mapsto \xi_z-\xi_{s}
\end{array}\right.\ \ \ \ \ \ \ \ \text{and}\ \ \ \ \ \ \ \ \left\{\begin{array}{l} \xi_4\mapsto \xi_t-\xi_{s}\\
\xi_5\mapsto \xi_u-\xi_{s}
\end{array}\right.\]
and
\[\left\{\begin{array}{l} \xi_1\mapsto \xi_x-\xi_s\\
\xi_2\mapsto \xi_x-\xi_{t}\\
\xi_3\mapsto \xi_x-\xi_{u}
\end{array}\right.\ \ \ \ \ \ \ \ \text{and}\ \ \ \ \ \ \ \ \left\{\begin{array}{l} \xi_4\mapsto \xi_x-\xi_{y}\\
\xi_5\mapsto \xi_x-\xi_{z}
\end{array}\right.\,.\]
This gives the $72$ transformations. The first set corresponds to the transformations sending $X_{\tm}$ to a matrix equivalent to itself, while the second set gives a minus sign. Permutations of $\{\xi_1,\xi_2,\xi_3\}$ and of $\{\xi_4,\xi_5\}$ are included in these $72$ transformations (choose $s=6$ in the first set).
\end{remax}

\subsection{Similarities with the magic square for the Clebsch--Gordan coefficients of $\mathfrak{su}_2$.\label{sec:CG}} 

In Section \ref{sec:CGsu2}, a connection between the Clebsch--Gordan coefficients of $\mathfrak{su}_2$ and the operators $X_\tm$ and $Y_\tm$ has already been pointed out.
Here, we establish the correspondence between the symmetry just identified and the well-known symmetry of the Clebsch--Gordan coefficients of $\mathfrak{su}_2$.
To make this statement precise, let us first recall the symmetry of the Clebsch--Gordan coefficients of $\mathfrak{su}_2$.

The $3jm$- (or Wigner) symbols are related the Clebsch--Gordan coefficients of $\mathfrak{su}_2$, as follows:
\begin{equation}
 \begin{pmatrix}
 J_1 &  J_2 & J_3\\
 M_1 & M_2 & -M_3
 \end{pmatrix}= \sqrt{2J_3+1} (-1)^{J_3+M_3}\langle J_1J_2; M_1  M_2 \ | \ J_1J_2; J_3 M_3 \rangle.
\end{equation}
Let us recall that $M_1+M_2=M_3$ for non-vanishing coefficients.
The symmetry of the Wigner coefficients can be read from the following magic square
\begin{equation}\label{eq:ms1}
 \left[\begin{array}{ccc}
2J_1+1 & 2J_2+1 &  2 J_3+1                        \\[0.5em]
J_2+J_3-M_1+1 & J_1+J_3-M_2+1 & J_1+J_2+M_3  +1          \\[0.5em]
J_2+J_3+M_1+1 & J_1+J_3+M_2+1 & J_1+J_2-M_3 +1
\end{array}\right].
\end{equation}
The sums of the entries in the same row or column are always equal to $2(J_1+J_2+J_3)+3$ (we used $M_1+M_2=M_3$).
The value of the Wigner coefficients is invariant under even permutations of the rows and columns, but also by transposition.
It is multiplied by $(-1)^{J_1+J_2+J_3}$ for odd permutations of the rows and columns.
Let us remark that we have added $-J_1-J_2-J_3-1$ to all the entries (and changed a global sign) in comparison to the usual magic square.
These transformations do not change the symmetries.

This magic square \eqref{eq:ms1} can be rewritten by using the parameters defined by \eqref{eq:eta} as follows 
\begin{equation}\label{eq:ms2}
 \left[\begin{array}{ccc}
\eta_1-\eta_5 & \eta_2 &  \eta_3-\eta_4                        \\[0.5em]
{}\eta_2-\eta_4 & \eta_3-\eta_5 & \eta_1       \\[0.5em]
{}\eta_3 & \eta_1-\eta_4 & \eta_2-\eta_5   
\end{array}\right].
\end{equation}
By comparing the magic square \eqref{eq:ms2} with one side of \eqref{square-xi}, the similarities between the symmetries of 
the Clebsch--Gordan coefficients of $\mathfrak{su}_2$ and those of the missing label operators $X$ and $Y$ become obvious. 

\begin{remax}
We do not have a satisfying explanation for this similarity between the symmetries of the missing label and those of the Clebsch--Gordan coefficients. 
Note for example that it is not obvious that the complicated formulas given in Proposition \ref{lem1} result in matrices $X_{\tm}$ and $Y_{\tm}$ having this symmetry.
\end{remax}

\section{The algebra $Z_2(\mathfrak{su}_3)_{\tm}$ and the symmetry of type $E_6$ \label{sec-alg}}

Recall from Section \ref{sec-XY} that $Z_2(\mathfrak{su}_3)$ is the centraliser in $U(\mathfrak{su}_3)^{\otimes 2}$ of the diagonal embedding of $\mathfrak{su}_3$. It is generated by $k_1,k_2,k_3,l_1,l_2,l_3,X,Y$. When the Casimir elements $k_1,k_2,k_3,l_1,l_2,l_3$ are replaced by numbers as in Remark \ref{rem-Cas}, the resulting algebra is denoted $Z_2(\mathfrak{su}_3)_{\tm}$ and is generated by $X$ and $Y$. 

In this section we will briefly recall the description of the algebra $Z_2(\mathfrak{su}_3)_{\tm}$ obtained in \cite{CPV1} in order to describe the relations satisfied by the tridiagonal matrices $X_{\tm}$ and $Y_{\tm}$ of Section \ref{sec-XY}. We shall give an interpretation of the symmetries of these matrices found in Theorem \ref{theo-sym} in terms of the Weyl group of type $E_6$ appearing in the description of $Z_2(\mathfrak{su}_3)_{\tm}$.

\subsection{Relations satisfied by $X_{\tm}$ and $Y_{\tm}$}

The matrices $X_{\tm}$ and $Y_{\tm}$ found in Section \ref{sec-XY} are the images of $X$ and $Y$ in the representation of the specialised 
centraliser $Z_2(\mathfrak{su}_3)_{\tm}$ on the multiplicity space $M_{\tm}$. They satisfy the following relations \cite{CPV1}:
\begin{subequations}\label{eq:def1}
\begin{align}
&[X_{\tm},Y_{\tm}]=Z_{\tm}\,, \\
&[X_{\tm},Z_{\tm}]= -6 Y_{\tm}^2 + a_2 X_{\tm}^2+ a_5 X_{\tm}+a_8\text{Id}\,,\\
&[Y_{\tm},Z_{\tm}]= -2X_{\tm}^3-a_2\{X_{\tm},Y_{\tm}\} -a_5 Y_{\tm}+a_6X_{\tm}+a_9\text{Id}\,,
\end{align}
\end{subequations}
where  $a_2,a_5,a_6,a_8,a_9,a_{12}\in\bC[\tm]$ (polynomials in the variables $m_1,m_2,m'_1,m'_2,m''_1,m''_2$).
The first relation defines the matrix $Z_{\tm}$. In addition to these commutation relations, we have: 
\begin{equation}
\label{eq:def2}
x_1X_{\tm}+x_2Y_{\tm}+x_3X_{\tm}^2+x_4\{X_{\tm},Y_{\tm}\} +x_5 Y_{\tm}^2+x_7 X_{\tm}Y_{\tm}X_{\tm}-X_{\tm}^4+4Y_{\tm}^3+Z_{\tm}^2=a_{12}\text{Id}\,,
\end{equation}
where the parameters $x_i$ are given by
\begin{equation*}
 x_1= 6 a_5  + 2 a_9\,,\ \ \ x_2= -2a_6 - 2 a_8\,,\ \ \ x_3= 6 a_2 + a_6\,,\ \ \ x_4=-a_5\,,\ \ \ x_5= 8 a_2 -24\,,\ \ \ x_7= -2 a_2 +12\ .
\end{equation*}
The commutation relations (\ref{eq:def1}) together with the special relation (\ref{eq:def2}) form a complete set of defining relations for the specialised centraliser 
$Z_2(\mathfrak{su}_3)_{\tm}$. We give in the next section an explicit expression for the coefficients $a_2,a_5,a_6,a_8,a_9,a_{12}$ as invariant polynomials for the action of a certain Weyl group of type $E_6$. The indices of these coefficients indicate their degrees as polynomials in the parameters in $\tm$.

\begin{remax}
The commutation relations (\ref{eq:def1}) are similar to the relations defining the Racah algebra, but they involve more general expressions in the right hand sides (in particular, a cubic polynomial). The last relation (\ref{eq:def2}) is the analogue of the relation that one must add to the Racah algebra to obtain the centraliser of the diagonal $\mathfrak{su}_2$ in $U(\mathfrak{su}_2)^{\otimes 3}$ (see \emph{e.g.} \cite{CFGPRV} or \cite{CGPV} for  the higher Racah algebras). Note that the left hand side of (\ref{eq:def2}) is central in the algebra defined by Relations (\ref{eq:def1}).
\end{remax}

\begin{remax} 
The matrices $H_1$ and $H_2$ defined in Section \ref{sec:CGsu2} satisfy also some relations, defining the so-called Hahn algebra, similar to the ones of (\ref{eq:def1}).
More precisely, $H_1$ and $H_2$ satisfy
\begin{eqnarray}
&&[H_1,H_2]=H_3\ ,\qquad [H_3,H_2]= 2 \{H_1,H_2\} +A_3 \ ,\qquad  [H_1,H_3]= 2 H_1^2 + H_2 +A_2\ ,
\end{eqnarray}
where the structure constants of the algebra are expressed in terms of $\eta$ as follows
\begin{subequations}\label{eq:AA}
\begin{eqnarray}
&& A_2=\frac{1}{2}(1-\eta_1^2-\eta_2^2-\eta_5^2-\eta_6^2)+\frac{1}{8}(\eta_1+\eta_2+\eta_5+\eta_6)^2 \ , \\
&&  A_3=-\frac{1}{4} (\eta_1+\eta_2-\eta_5-\eta_6) \Big( (\eta_1-\eta_2)^2-(\eta_5-\eta_6)^2 \Big) \ .
\end{eqnarray}
\end{subequations}
\end{remax}

\subsection{Invariance under the Weyl group of type $E_6$}  \label{sec:E6}

Let us consider the root system of type $E_6$. We fix a set of simple roots $\alpha_1,\alpha_2,\alpha_3,\alpha_4,\alpha_5,\alpha_6$ in the 
Euclidean vector space $\mathbb{R}^6$ corresponding to the following Dynkin diagram and the following Cartan matrix $A_{ij}=(\alpha_i,\alpha_j)$:
\begin{center}
\begin{tikzpicture}[scale=0.2]
\draw (1,1) circle [radius=1];
\node [above] at (1,1.5) {$1$};
\draw (2,1)--(4,1);

\draw (5,1) circle [radius=1];
\node [above] at (5,1.5) {$2$};
\draw (6,1)--(8,1);

\draw (9,1) circle [radius=1];
\node [above] at (9,1.5) {$3$};
\draw (10,1)--(12,1);

\draw (13,1) circle [radius=1];
\node [above] at (13,1.5) {$4$};
\draw (14,1)--(16,1);

\draw (17,1) circle [radius=1];
\node [above] at (17,1.5) {$5$};

\draw (9,-3) circle [radius=1];
\node [right] at (9.5,-3) {$6$};
\draw (9,0)--(9,-2);
\node [right] at (27,1) {$A=\begin{pmatrix}
    2& -1 &0 &0 &0 &0\\
    -1 & 2 & -1 & 0 & 0 & 0 \\
    0&-1&2&-1&0&-1\\
    0&0&-1&2&-1&0\\
    0&0&0&-1&2&0\\
    0&0&-1&0&0&2
   \end{pmatrix}\ .$};
\end{tikzpicture}
\end{center}
The Weyl group $W(E_6)$ is generated by the simple reflections $\{s_i\}_{i=1,\dots,6}$, which act on $\mathbb{R}^6$ by:
\begin{equation}\label{act-W}
s_i(x)=x-(x,\alpha_i)\alpha_i\ \ \ \ \ \forall x\in\mathbb{R}^6\ .
\end{equation}
Let us associate the parameters in $\tm$ with some roots as follows:
\begin{equation}\label{eq:relrp}
 \begin{array}{lll}
 m_1=\alpha_3+\alpha_4+\alpha_5,\ \  &\ \  m'_1=\alpha_1+\alpha_2+\alpha_3, \ \  &\ \  m''_1=\alpha_3+\alpha_6, \\[0.4em]
m_2=\alpha_1+\alpha_2+\alpha_3+\alpha_4+\alpha_6, \ \  &\ \  m'_2=\alpha_2+\alpha_3+\alpha_4+\alpha_5+\alpha_6, \ \  &\ \  m''_2=\alpha_2+\alpha_3+\alpha_4,
\end{array}
\end{equation}
\begin{remax}
One can act with an element of the Weyl group $W(E_6)$ on the right-hand-sides in (\ref{eq:relrp}) and obtain another correspondence which would not change the expression of the coefficients below. 
The choice we make here will turn out to be well adapted later for the description of symmetries. Up to the action of $W(E_6)$, 
it is equivalent to the canonical choice made in \cite{CPV1}: $(m_1,m_2,m'_1,m'_2,m''_1,m''_2)\mapsto (\alpha_1,\alpha_1,\alpha_5,\alpha_4,\Theta,-\alpha_6)$ where 
$\Theta=\alpha_1+2\alpha_2+3\alpha_3+2\alpha_4+\alpha_5+2\alpha_6$ is the longest positive root.
\end{remax}

The Weyl group $W(E_6)$ acts on the set of roots via (\ref{act-W}) and thus acts also linearly on the parameters $m_1,\dots,m''_2$ through the correspondence above. 
Using this action, we introduce the homogeneous polynomials $p_i\in\bC[\tm]$ of degree $i$ given by:
\begin{subequations}
\begin{eqnarray}
&& p_2=\frac{3}{2} \langle m_1^2 \rangle\ ,\quad p_5=\frac{8}{3}\langle (m'_1)^2m'_2m''_1m''_2 \rangle\ ,\quad p_6=10\langle m_1m_2m'_1m'_2m''_1m''_2\rangle\,,\\
&&  p_8=\frac{5}{3} \langle (m_1)^2m_2(m'_1)^2m'_2m''_1m''_2 \rangle \ ,\quad p_9=\frac{40}{27}\langle (m_1m_2m'_1)^2m'_2m''_1m''_2 \rangle\,,\\
&&  p_{12}=\frac{20}{3}\langle( m_1m_2m'_1m'_2m''_1m''_2)^2\rangle\,,
\end{eqnarray}
\end{subequations}
where the average $\langle.\rangle$ is defined as follows
\begin{equation}
 \langle. \rangle=\frac{1}{51840} \sum_{s \in W(E_6)} s(.) \ .
\end{equation}
The normalisation factor $51840$ is the order of the Weyl group $W(E_6)$. 
These polynomials are independent and their degrees $2,5,6,8,9,12$ are the fundamental degrees of type $E_6$. 
Thus, they form a set of homogeneous generators for the invariant polynomials of the Weyl group $W(E_6)$.

Finally, it is found in \cite{CPV1} that the parameters appearing in  (\ref{eq:def1}) and (\ref{eq:def2})  are expressed in terms of these invariant polynomials according to:
\begin{subequations}
\begin{eqnarray}
  && a_2=p_2 -3 \ , \quad a_5 =-p_5 \ , \quad a_6=p_6+  \frac{p_2^3}{9}+\frac{2p_2^2}{3}-\frac{3p_2}{2}+1\,,\\
  && a_8=-p_8 +\frac{p_2^4}{54} +\frac{p_2p_6}{12}  +\frac{p_2^3}{18}+\frac{p_6}{2}+\frac{p_2^2}{6}-\frac{p_2}{4}+\frac{1}{8}\,,\\
  &&a_9 = -p_9-p_5\left(\frac{p_2^2}{27}+\frac{p_2}{3}- \frac{1}{4}\right)\,,
 \end{eqnarray}
and
\begin{eqnarray}
 a_{12}&=& -p_{12}+\frac{35p_6^2}{12} +\frac{p_2^6}{36}+\frac{17p_2^3p_6}{72} -\frac{p_2^2p_8 }{18} -\frac{7p_2p_5^2}{18}
  +\frac{p_2^5}{162}-\frac{p_2p_8}{3}+\frac{p_2^2p_6}{36}-\frac{p_5^2}{4}\\
  &&\hspace{1cm}-\frac{13p_2^4}{108}+\frac{13p_8}{2}-\frac{13p_2p_6}{24}-\frac{19p_2^3}{54}-3p_6-\frac{11p_2^2}{12}+\frac{11p_2}{8}-\frac{11}{16}\,.\nonumber
\end{eqnarray}
\end{subequations}
As a consequence of this invariance of the coefficients, it is immediate that the relations in (\ref{eq:def1}) and (\ref{eq:def2}) are preserved by 
the group $W(E_6)$. In fact, if we denote $\pm W(E_6)$ the group generated by $W(E_6)$ and the central symmetry $-\text{Id}_{\mathbb{R}^6}$, we have the following theorem.
\begin{theo} \label{th:symE6}
If $\tm'=s.\tm$ with $s\in\pm W(E_6)$ then  $Z_2(\mathfrak{su}_3)_{\tm}\cong Z_2(\mathfrak{su}_3)_{\tm'}$.
The isomorphism is given by $X\mapsto \pm X$ and $Y\mapsto Y$ (the sign corresponds to the sign of $s\in\pm W(E_6)$).
\end{theo}

\begin{remax}
The centraliser $Z_2(\mathfrak{su}_3)_{\tm}$ is studied with a different approach through a certain quantum Hamiltonian reduction in \cite{ELOR}. In particular, $Z_2(\mathfrak{su}_3)_{\tm}$ is related to the spherical subalgebra of a certain symplectic reflection algebra of rank 1.
\end{remax}

\subsection{Symmetry of the missing label inside $\pm W(E_6)$}\label{subsec-symE6}

Through the correspondence (\ref{eq:relrp}) between parameters and roots, the set of $18$ elements in (\ref{eq:deg}) of which the minimum gives the multiplicity $d_{\tm}$ becomes:
\begin{equation}\label{eq:deg-roots}\left\{\begin{array}{cccccc}
345\,, & 123\,, & 234\,,              \quad              &\quad 3456\,, & 1236\,, & 2346\,, \!\!  \\[0.5em]
\!\!23\,,\!\! & 34\,,\!\! & 12345\,,  \quad    & \quad   236\,,\!\! & 346\,,\!\! & 123456\,, \\[0.5em]
1234\,, & 2345\,, & 3\,,  \quad   & \quad 12346\,, & 23456\,, & 36
\end{array}\right\}\ ,
\end{equation} 
where we use the short-hand notation $ij...k$ for the root $\alpha_i+\alpha_j+\dots+\alpha_k$. 
The remarkable fact, justifying our choice of correspondence in (\ref{eq:relrp}), is that these 18 roots are exactly all the positive roots of $E_6$ 
which have $\alpha_3$ with coefficient 1, see Appendix \ref{app-diagram}.

Our goal in this section is to identify the subgroup of the Weyl group $\pm W(E_6)$ leaving this set of roots stable. 
This gives another description of the symmetries of the missing label in Theorem \ref{theo-sym}. 
This subgroup stabilizing (\ref{eq:deg-roots}) can be computed by testing all the elements of $\pm W(E_6)$ but 
this is a lengthy check and an alternative way is given below. 

Let $(\omega_1,\dots,\omega_6)$ denote the basis of $\textbf{R}^6$ dual to the basis of simple roots. 
Then the set of roots above is the set of roots $\beta$ such that $(\beta,\omega_3)=1$. 
So an element of $W(E_6)$ which stabilises $\omega_3$ will leave the above subset stable. 
Reciprocally, the elements $\pi$ of $W(E_6)$ leaving the above set stable must satisfy $\bigl(\beta,\omega_3-\pi(\omega_3)\bigr)=0$ for any root $\beta$ in the above set. 
This implies that $\pi(\omega_3)=\omega_3$. We thus conclude that the subgroup of $W(E_6)$ which leaves the above set of roots stable coincides with the stabiliser of $\omega_3$. 
This stabiliser is a parabolic subgroup of $W(E_6)$ and it is well-known that it is generated by the reflections it contains. So it is easy to see that it is the subgroup of $W(E_6)$ generated by 
all simple reflections except $s_3$.  

We must also add the central symmetry $-\text{Id}$ into the picture. Let $w_0$ denote the longest element of $W(E_6)$.
We recall that $w_0$ is of order $2$ and that the element $-w_0$ corresponds to the automorphism of the Dynkin diagram of $E_6$ that 
exchanges $\alpha_1,\alpha_2$ with, respectively, $\alpha_5,\alpha_4$ (and leaves $\alpha_3$ and $\alpha_6$ invariant). Thus, the element $-w_0$ also leaves the above set of roots stable.

We conclude that the subgroup of $\pm W(E_6)$ preserving the set (\ref{eq:deg-roots}) is generated by the elements $s_1,s_2,s_4,s_5,s_6,-w_0$. As recalled above, we have $w_0s_1w_0=s_5$ and $w_0s_2w_0=s_4$. Moreover, $s_6$ commutes with all other generators, while $s_1,s_2$ commute with $s_4,s_5$. So the structure of this group is:
\begin{equation}\label{group2}\Bigl((\langle s_1,s_2\rangle\times \langle s_4,s_5 \rangle)\rtimes \langle -w_0\rangle\Bigr)\times \langle s_6\rangle\ .
\end{equation}
The subgroups $\langle s_1,s_2\rangle$ and $\langle s_4,s_5 \rangle$ are isomorphic to symmetric groups $S_3$ and 
so this description is reminiscent of the symmetry group of the missing label obtained in (\ref{group1}).

Using the diagram of Appendix \ref{app-diagram}, it is straightforward to describe the action of the different generators on the set (\ref{eq:deg-roots}) 
(which as we recall is the array of the two magic squares given in the introduction when parameters are identified with roots according to (\ref{eq:relrp})):
\begin{itemize}
\item $s_1$ is the transposition and $s_2$ is the reflection through the diagonal $(123,12345,1234)$ (both simultaneously in the two $3\times 3$ square);
\item $s_4$ is the reflection through the antidiagonal $(345,12345,2345)$ and $s_5$ is the reflection through the antidiagonal $(123,23,3)$ (both simultaneously in the two $3\times 3$ square);
\item $-w_0$ is the transposition of the first two columns (simultaneously in the two $3\times 3$ square) and $s_6$ is the exchange of the two squares.
\end{itemize}
It is easy to see that the group generated by these transformations coincides with the group of transformations described in Theorem \ref{theo-sym}. 
This leads to the following alternative description of the symmetries of the missing label:
\begin{coro}\label{coro-sym}
Thanks to the identification (\ref{eq:relrp}) between parameters and roots of $E_6$, the group of symmetries of the missing label corresponds to the subgroup of $\pm W(E_6)$ generated by:
\[s_1,s_2,s_4,s_5,s_6\ \text{and}\ -w_0\ .\]
For a transformation $g\ :\ \tm\mapsto \tm'$ in this subgroup, the pair of matrices $(X_{\tm},Y_{\tm})$ is conjugated to $(\pm X_{\tm'},Y_{\tm'})$, the sign being $(-1)$ if and only if $g\in -W(E_6)$.
\end{coro}

\section{Representations of the missing label algebra \label{sec:rep}}

In this section, we show that, by using the root system $E_6$, there exists a geometric way to construct 
some representations of the algebra $Z_2(\mathfrak{su}_3)_{\tm}$ satisfied by the missing label operators.
We show that to each of the $36$ positive roots of $E_6$, we can associate $12$ representations involving a pair of infinite-dimensional tridiagonal matrices. Each of these representations is associated to a different $4$-simplex which is a $4$-face of the root polytope of $E_6$. We also show how to extract finite-dimensional representations of $Z_2(\mathfrak{su}_3)_{\tm}$ from these infinite-dimensional representations. The physical representation of Section \ref{sec-XY} is identified as one of these.

\subsection{Infinite-dimensional representations and $4$-faces}\label{subsec-4face}

In this subsection, we consider the general situation where the parameters in $\tm$ consist of $6$ independent complex numbers
and $Z_2(\mathfrak{su}_3)_{\tm}$ is the algebra with generators $X,Y$ and relations as in (\ref{eq:def1})-(\ref{eq:def2}), where the subscript $\tm$ is dropped. 
We keep the identification (\ref{eq:relrp}) between the 6 parameters $\tm$  and the roots of $E_6$.

\paragraph{4-faces of the root polytopes.} Let $\Theta=\alpha_1+2\alpha_2+3\alpha_3+2\alpha_4+\alpha_5+2\alpha_6$  be the highest root of $E_6$
and let us associate to this root the (unordered) set $\Xi_{\Theta,0}$ of $5$ other positive roots $\{\xi^0_1,\xi^0_2,\xi^0_3,\xi^0_4,\xi^0_5\}$ given by:
\[\Xi_{\Theta,0}=\{\alpha_1+\alpha_2+\alpha_3+\alpha_4+\alpha_5,\ \alpha_1+\alpha_2+\alpha_3+\alpha_4,\ \alpha_1+\alpha_2+\alpha_3,\ \alpha_1+\alpha_2,\ \alpha_1\}\ .\]
For any positive root $\Lambda$, there exists an element $w\in W(E_6)$ of minimal length such that $w(\Theta)=\Lambda$ which permits to define 
the following set:
\[\Xi_{\Lambda,0}=\{\xi_1,\xi_2,\xi_3,\xi_4,\xi_5\}=\{w(\xi^0_1),\,w(\xi^0_2),\,w(\xi^0_3),\,w(\xi^0_4),\,w(\xi^0_5)\}\ .\]
It is straightforward to calculate them using the diagram for the root poset given in Appendix \ref{app-diagram}. Indeed the element $w$ corresponds to a shortest path from $\Theta$ to $\Lambda$ in the diagram, and the action of an element $w$ on any root is easily read off the diagram.
Two elements $w,w'$ satisfy $w(\Theta)=w'(\Theta)$ if and only if they belong to the same coset, 
$wG_{\Theta}=w'G_{\Theta}$, where $G_{\Theta}$ is the stabiliser subgroup of $\Theta$. 
With the same kind of reasoning as in Section \ref{subsec-symE6}, we find that $G_{\Theta}$ is the parabolic subgroup generated by $s_1,s_2,s_3,s_4,s_5$. 
Therefore in each coset, there is a unique element $w$ of minimal length, which moreover satisfies $\ell(w\pi)=\ell(w)+\ell(\pi)$ for all $\pi\in G_{\Theta}$ \cite[Proposition 2.1.1]{GP}. 
This shows the unicity of two elements $w$ and, in addition, we have that $\ell(ws_i)=\ell(w)+1$ for all $i=1,\dots,5$. This is well known to be equivalent to the fact that $w$ sends $\alpha_i$ to a positive root 
\cite[Theorem 1.1.9]{GP}. We conclude that $\xi_1,\xi_2,\xi_3,\xi_4,\xi_5$ as defined above are positive roots.

\begin{exam}\label{ex:4fa} To offer some examples, we refer to the diagram in Appendix \ref{app-diagram}. 
First, let $\Lambda=\Theta-\alpha_6$. Then, the element $w$ is $s_6$ and we find that $\Xi_{\Theta-\alpha_6,0}$ is:
\[\{\alpha_1+\alpha_2+\alpha_3+\alpha_4+\alpha_5+\alpha_6,\ \alpha_1+\alpha_2+\alpha_3+\alpha_4+\alpha_6,\ \alpha_1+\alpha_2+\alpha_3+\alpha_6,\ \alpha_1+\alpha_2+\alpha_6,\ \alpha_1+\alpha_6\}\ .\]
Second, let $\Lambda=\alpha_1$. Then, the element $w$ is $s_2s_3s_4s_5s_6s_3s_2s_4s_3s_6$ and we find:
\[ \Xi_{\alpha_1,0}=\{ \Theta,\ \Theta-\alpha_6,\ \Theta-\alpha_6-\alpha_3,\ \Theta-\alpha_6-\alpha_3-\alpha_4,\ \Theta-\alpha_6-\alpha_3-\alpha_4-\alpha_5 \}\ .\]
\end{exam}

Let us make some remarks about the sets $\Xi_{\Lambda,0}$:
\begin{itemize}
 \item the five roots in $\Xi_{\Lambda,0}$ form a $4$-simplex which is a $4$-face of the polytope formed by the roots in $E_6$ (there are 432 such $4$-faces in total, see \cite{Sza}). We have constructed so far 36 such $4$-faces.
 \item It is easily seen that the $4$-face $\Xi_{\Theta,0}$ is orthogonal to $\Theta$, and moreover, that $\Theta$ is the unique positive root orthogonal to $\Xi_{\Theta,0}$. Since the elements of $W(E_6)$ preserve the scalar product, we have that the $4$-face $\Xi_{\Lambda,0}$ is orthogonal to a unique positive root, which is $\Lambda$.
\item In a given $\Xi_{\Lambda,0}$, the roots $\xi_i$ satisfy $(\xi_i,\xi_j)=1$ for $1\leq i\neq j \leq 5$. 
Indeed it is true for the particular $4$-face $\Xi_{\Theta,0}$, and remains true for any $\Lambda$ since the elements of $W(E_6)$ preserve the scalar product.
\end{itemize}

To construct all the $432$ $4$-faces which are $4$-simplices, we proceed as follows. 
We start with one of the $4$-faces $\Xi_{\Lambda,0}=\{\xi_1,\xi_2,\xi_3,\xi_4,\xi_5\}$ as defined above, 
we choose some $k\in\{1,\dots,5\}$ and we apply the reflection corresponding to the root $\xi_k$ on each element of $\Xi_{\Lambda,0}$. This yields:
\begin{equation}
 \Xi_{\Lambda,k}=\{\ -\xi_k,\quad \xi_j-\xi_k\ | \ j\neq k\ \}.
\end{equation}
Note that the six $4$-faces $\Xi_{\Lambda,k}$ (for $0\leq k\leq 5$) obtained this way are pairwise different and are orthogonal to the same positive root $\Lambda$. At this point, 
we have constructed $216$ $4$-faces which form a single orbit under the action of the Weyl group $W(E_6)$. The remaining half is provided by the images under the central symmetry $-\text{Id}$. 
Namely is is made out of the $4$-faces
\begin{equation}
 -\Xi_{\Lambda,k}=\{\ -\xi \ |\ \xi \in \Xi_{\Lambda,k}\ \}.
\end{equation}
As is obvious from the construction, the $432$ $4$-faces form a single orbit under the group $\pm W(E_6)$. 

\paragraph{Construction of representations.}
Take a positive root $\Lambda$ and a 4-face $\mathcal{F}=\{\xi_1,\xi_2,\xi_3,\xi_4,\xi_5\}$ orthogonal to $\Lambda$ as explained above. We also need to introduce $\lambda_+$ and $\lambda_-$ defined by 
\[\lambda_{\pm}=\pm\frac{1}{2}\Lambda+\frac{1}{6}(\xi_1+\xi_2+\xi_3+\xi_4+\xi_5)\ .\]
We will use the notations $\xi_6=0$ and
\[E_{j,k}=\sum_{i=1}^6x_i(j)^k\,,\ \ \ \ \text{where $x_i(j)=(\xi_i-j+\frac{1}{2})$\ for $i\in\{1,\dots,6\}$\ .}\]
We are now in a position to provide representations of the algebra $Z_2(\mathfrak{su}_3)_{\tm}$. 
We define tridiagonal matrices $\cA^{\mathcal{F}}$ and $\cB^{\mathcal{F}}$, whose rows and columns are indexed by $j\in\mathbb{Z}$, and whose non-vanishing entries are given by:
\begin{subequations}
\begin{eqnarray}
&& \cA^{\mathcal{F}}_{j+1,j}= 1 \ , \qquad\qquad \cA^{\mathcal{F}}_{j,j+1}=j(j-\xi_1)(j-\xi_2)(j-\xi_3)(j-\xi_4)(j-\xi_5)  \,, \label{eq:cX}\\
&& \cB^{\mathcal{F}}_{j+1,j}=(j-\lambda_+)\,,  \quad  \cB^{\mathcal{F}}_{j,j+1}= j(j-\xi_1)(j-\xi_2)(j-\xi_3)(j-\xi_4)(j-\xi_5)(j-\lambda_-)\label{eq:cY}\ ,
\end{eqnarray}
\end{subequations}
and 
\[ \cA^{\mathcal{F}}_{jj}=-\frac{1}{108}\bigl(\frac{7}{2}E_{j,1}^3-18 E_{j,1}E_{j,2}+18 E_{j,3}\bigr)-\frac{1}{24}E_{j,1}\bigl(\Lambda^2+2\bigr)\ ,\]
\[\cB^{\mathcal{F}}_{jj}=\frac{1}{288}\Bigl(\frac{5}{2}E_{j,1}^4+32 E_{j,1}E_{j,3}+6\bigl(E_{j,2}^2-3E_{j,1}^2E_{j,2}-4E_{j,4})+6E_{j,2}(\Lambda^2+2)-3E_{j,1}^2(\Lambda^2-2)-\frac{3}{2}\Lambda^4+6\Lambda^2-36\Bigr)\,.\]
\begin{thm}\label{th:pr}
Fix a positive root $\Lambda$, an integer $k\in\{0,\dots,5\}$ and a sign $\epsilon\in\{-1,1\}$, 
and let $\mathcal{F}=\epsilon\Xi_{\Lambda,k}$ be the associated $4$-face. Set $X^{\mathcal{F}}=\epsilon\cA^{\mathcal{F}}$ and $Y^{\mathcal{F}}=\cB^{\mathcal{F}}$.\\
The assignment $X\mapsto X^{\mathcal{F}}$ and $Y\mapsto Y^{\mathcal{F}}$ provides a representation of the 
algebra $Z_2(\mathfrak{su}_3)_{\tm}$, that we denote $V_{\mathcal{F}}^{\infty}$.
\end{thm}
\proof Assume that the theorem is satisfied for some $4$-face $\mathcal{F}$ and take some element $s\in\pm W(E_6)$. 
The matrices $\cA^{s\cdot\mathcal{F}}$ and $\cB^{s\cdot\mathcal{F}}$ satisfy the relations of 
the algebra $Z_2(\mathfrak{su}_3)_{s\cdot\tm}$. Recalling that $Z_2(\mathfrak{su}_3)_{s\cdot\tm}$ is isomorphic to $Z_2(\mathfrak{su}_3)_{\tm}$ 
(see Theorem \ref{th:symE6}), this proves the statement of the theorem for the $4$-face $s\cdot\mathcal{F}$. 
Therefore, since the set of $4$-faces forms a single orbit under the action of $\pm W(E_6)$, it is enough to check the assertion 
of the theorem for a single $4$-face (for example, the $4$-face $\Xi_{\Theta,0}$). The verification of 
the defining relations \eqref{eq:def1}-\eqref{eq:def2} of $Z_2(\mathfrak{su}_3)_{\tm}$ is straightforward and 
can be done as follows. If we denote by $v_j$ the basis vector indexed by an integer $j\in\mathbb{Z}$, then the defining relations 
\eqref{eq:def1}-\eqref{eq:def2} applied on these $v_j$ only involves the indices between $j-4$ and $j+4$. It thus suffices to 
consider the $9$ by $9$ matrices formed by the entries with indices in that range and to verify directly the defining relations. This calculation can be made in a straightforward manner by hand (very lengthy) or immediately using a computer.
\endproof

\begin{remax}\label{rem-trans}
In the formulas for $\cB^{\mathcal{F}}$ given above, a choice was made regarding the role of $\lambda_+$ and $\lambda_-$.
It turns out that exchanging $\lambda_+$ and $\lambda_-$ also leads to representations of $Z_2(\mathfrak{su}_3)_{\tm}$. 
One can show that these other representations are equivalent to the transpose of the representations that we already have. 
In fact it is straightforward to check that exchanging $\lambda_+$ and $\lambda_-$ in $V_{\mathcal{F}}^{\infty}$ is equivalent to transposing the matrices in $V_{-\mathcal{F}}^{\infty}$, up to the renaming of indices $j\mapsto -j+1$.
\end{remax}

\subsection{Finite-dimensional representations}

The representations constructed in the previous subsection are infinite-dimensional but,
thanks to the tridiagonal property of the representing matrices, finite-dimensional representations can be extracted. 

For simplicity, we assume from now on that all parameters are integers:
\[\tm\in(\mathbb{Z})^6\ \ \ \ \ \text{and}\ \ \ \ \ n,\ell\in\mathbb{Z}\ .\]
These conditions are satisfied, in particular, if the parameters $\tm$ are positive integers: $\tm\in(\mathbb{Z}_{>0})^6$ such that the irreducible representation $[m''_1,m''_2]$ appears in the tensor product $[m_1,m_2]\otimes[m'_1,m'_2]$ with non-zero multiplicity (the situation of the preceding Sections).

From the association between parameters and roots of $E_6$ given in (\ref{eq:relrp}), we see that to every roots corresponds a value in $\mathbb{Z}$. 
Therefore for any $4$-face $\mathcal{F}=\{\xi_1,\dots,\xi_5\}$, the parameters $\xi_1,\dots,\xi_5$ take values in $\mathbb{Z}$.

Let $\mathcal{F}=\{\xi_1,\xi_2,\xi_3,\xi_4,\xi_5\}$ be a $4$-face and denote (as in Theorem \ref{th:pr}) by $X^{\mathcal{F}}$ and $Y^{\mathcal{F}}$ 
the infinite matrices representing $X$ and $Y$ in the corresponding representation of $Z_2(\mathfrak{su}_3)_{\tm}$.

Recall that we have set $\xi_6:=0$ and now rename the parameters $\xi$'s so that:
 \[\xi_1\leq\xi_2\leq\xi_3\leq\xi_4\leq\xi_5\leq\xi_6\ .\]
Since they are integers, we can extract finite-dimensional representations as in the following definition.
\begin{defi}\label{def-rep}
Let $a\in\{1,...,5\}$. Set:
\[X^{\mathcal{F},a}=|X^{\mathcal{F}}_{i,j}|_{\xi_a+1\leq i,j \leq \xi_{a+1}}\ \ \ \text{and}\ \ \ Y^{\mathcal{F},a}=|Y^{\mathcal{F}}_{i,j}|_{\xi_a+1\leq i,j \leq \xi_{a+1}}\ .\]
$V^{(a)}_{\mathcal{F}}$ denotes the finite-dimensional representation of $Z_2(\mathfrak{su}_3)_{\tm}$ afforded by these two matrices: $X\mapsto X^{\mathcal{F},a}$ and $Y\mapsto Y^{\mathcal{F},a}$.
\end{defi}
This indeed defines representations of $Z_2(\mathfrak{su}_3)_{\tm}$ by the following argument. Since $X^{\mathcal{F}}_{\xi_a,\xi_{a}+1}=Y^{\mathcal{F}}_{\xi_a,\xi_a+1}=0$, 
the subspace spanned by the basis vectors with indices strictly superior to $\xi_a$ is stable. Similarly, the subspace spanned by the basis vectors with indices strictly superior to $ \xi_{a+1}$ 
is stable. We take the quotient of the first one by the second one. The resulting quotient representation is afforded by the matrices given in definition \ref{def-rep}.

\paragraph{Equivalent representations.} Consider the two $4$-faces $\mathcal{F}=\pm\Xi_{\lambda,0}$ and $\mathcal{F}'=\pm\Xi_{\lambda,k}$. Passing from $\mathcal{F}$ to $\mathcal{F}'$, 
the parameters $\{\xi_i\}_{i=1,\dots,6}$ are replaced by $\{\xi_i-\xi_k\}_{i=1,\dots,6}$ (recall that $\xi_6:=0$). 
Moreover, the parameter $\Lambda$ is invariant and the parameters $\lambda_\pm$ become $\lambda_\pm-\xi_k$. 
Since $\xi_k$ is an integer, the indices of the matrices $X^{\mathcal{F}'}$ and $Y^{\mathcal{F}'}$ can be shifted as follows $j\mapsto j+\xi_k$ to recover the matrices $X^{\mathcal{F}}$ and $Y^{\mathcal{F}}$.

We have just shown that it is enough to consider the infinite-dimensional representations associated to the $72$ following $4$-faces:
\[V^{\infty}_{\mathcal{F}}\,,\ \ \quad\text{where $\mathcal{F}\in\{\pm \Xi_{\Lambda,0}\}_{\text{$\Lambda$ a positive root}}$}.\] 
Note then from the defining relations (\ref{eq:def1})-(\ref{eq:def2}) of the algebra $Z_2(\mathfrak{su}_3)_{\tm}$ that if $X\mapsto \cA$ and $Y\mapsto \cB$ provide a matrix representation,
the transpose matrices ${}^t\cA$ and ${}^t\cB$ give a representation as well. As explained in Remark \ref{rem-trans}, one can check that the infinite-dimensional representation $V^{\infty}_{-\mathcal{F}}$ is equivalent, up to transposition, to $V^{\infty}_{\mathcal{F}}$, with the role of $\lambda_+$ and $\lambda_-$ exchanged. Therefore the matrices representing $X$ and $Y$ in a finite-dimensional representation extracted from $V^{\infty}_{-\mathcal{F}}$ are equivalent to matrices of the form:
\[\left(\begin{array}{cccc} x_{1,1} & 1 & & 0\\
x_{1,2} & x_{2,2} & \ddots & \\
 & \ddots & \ddots & 1\\
0 & & x_{d-1,d} & x_{d,d}\end{array}\right)\ \ \ \text{and}\ \ \ \left(\begin{array}{cccc}y_{1,1} & \tilde{y}_{2,1} & & 0\\
\tilde{y}_{1,2} & y_{2,2} & \ddots & \\
 & \ddots & \ddots & \tilde{y}_{d,d-1}\\
0 & & \tilde{y}_{d-1,d} & y_{d,d}\end{array}\right)\,,\]
where $(x_{i,j})$ and $(y_{i,j})$ are matrices extracted from $V^{\infty}_{\mathcal{F}}$ and $\tilde{.}$ means exchanging $\lambda_+$ and $\lambda_-$.

Now note that all coefficients $x_{j,j+1}$ are different from $0$ by construction of the extracted representations and the following equalities result from the explicit values of the off-diagonal coefficients:
\[\tilde{y}_{j,j+1}=x_{j,j+1}y_{j+1,j}\ \ \ \text{and}\ \ \ \tilde{y}_{j+1,j}=y_{j+1,j}/x_{j,j+1}\ .\]
We deduce easily that conjugating by a diagonal matrix $(d_{ij})$ with coefficients satisfying $d_{j+1,j+1}/d_{j,j}=x_{j,j+1}$, 
leads to the original matrices $(x_{i,j})$ and $(y_{i,j})$ that were obtained from $V^{\infty}_{\mathcal{F}}$.

We conclude that it is enough to consider the representations extracted from one of the $4$-faces $\Xi_{\Lambda,0}$. 
Namely, up to transposition, the finite-dimensional representations we have constructed are all equivalent to one of the following:
\[V^{(a)}_{\mathcal{F}}\,,\ \ \quad\text{where $a\in\{1,\dots,5\}$ and $\mathcal{F}\in\{\Xi_{\Lambda,0}\}_{\text{$\Lambda$ a positive root}}$}.\]

\paragraph{Identification of the physical representations.} 
Assume that $\tm\in(\mathbb{Z}_{>0})^6$ and that the multiplicity $d_{\tm}$ is non-zero. Consider first the case $n\leq\ell$. The correspondence (\ref{eq:relrp}) between parameters and roots is such that the set $\{\xi_1,\xi_2,\xi_3,\xi_4,\xi_5\}$ used in Section \ref{subsec-formulas} corresponds, in terms of roots, to the $4$-face:
\[\Xi_{\Theta,0}=\{\alpha_1+\alpha_2+\alpha_3+\alpha_4+\alpha_5,\ \alpha_1+\alpha_2+\alpha_3+\alpha_4,\ \alpha_1+\alpha_2+\alpha_3,\ \alpha_1+\alpha_2,\ \alpha_1\}\ ,\]
orthogonal to the highest root $\Theta$. Since $d_{\tm}>0$, we have:
\[\{\xi_4,\xi_5,0\}<\{\xi_1,\xi_2,\xi_3\}\ ,\]
and so we can extract matrices by taking indices from $\text{max}\{\xi_4,\xi_5,0\}+1$ to $\text{min}\{\xi_1,\xi_2,\xi_3\}$. This gives the matrices for the missing label operators $X_{\tm}$ and $Y_{\tm}$ found in Section \ref{subsec-formulas}.

Now if $\ell\leq n$, still using the correspondence (\ref{eq:relrp}) between parameters and roots, we need to consider instead the $4$-face orthogonal to the root $\Theta-\alpha_6$, for which we have found in Example \ref{ex:4fa}, $\Xi_{\Theta-\alpha_6,0}=\Xi_{\Theta,0}+\alpha_6$ (note that $\alpha_6$ corresponds to $\ell-n$).

In terms of $4$-faces, one can check that the 144 symmetries of the missing label operators correspond to transformations (say $n\leq \ell$) 
which either permute the roots in the $4$-face $\Xi_{\Theta,0}$, or which transforms this $4$-face into $\Xi_{\Theta,k}$ or $-\Xi_{\Theta,k}$. 
This extends to general statements about equivalent representations associated to any $4$-face.

 \section{Conclusion}
 
 We tackled the so-called missing label problem for $\mathfrak{su}_3$ by studying its diagonal centraliser generated by two operators $X$ and $Y$. 
 We obtained a tridiagonal representation of these operators whose eigenvalues provide the missing label. We studied extensively their symmetries,
 the associated integrable systems, showed how the Bethe ansatz applies and constructed explicit examples.
We further examined on detail the centralizer, its symmetry, provided by the Weyl group of $E_6$ and some of its representations labeled by the $4$-faces of the $E_6$-polytope.
 
 The interplay between centraliser, Weyl group, algebraic Heun operator and Bethe ansatz is fascinating and deserves further investigations.
 Firstly, the correlation functions of the Gaudin model should provide a way to compute the Clebsch--Gordan coefficients of $\mathfrak{su}_3$.
 Secondly, the study of the centraliser of the diagonal action for $U_q(\mathfrak{su}_3)$ in two fold tensor product certainly deserves attention. 
 We may obtain a $q$-deformation of the 
 algebra $Z_\tm(\mathfrak{su}_3)$ and wonder if the $W(E_6)$-symmetry is preserved and if the Bethe ansatz works also in this case.
 Thirdly, the description of the centraliser for higher rank Lie algebra is another open problem and should give some insights to identify the missing label(s).
 Finally, we would like to examine the centraliser of the diagonal action for $U_q(\mathfrak{su}_3)$ in the $n$ fold tensor product.
 We want to give a description in terms of generators and relations of this centraliser, to identify how the $W(E_6)$-symmetry is modified, 
 to find the missing labels and ti diagonalise them through the Bethe ansatz approach.
 \\
 
 \vspace{1cm}
\noindent
\textbf{Acknowledgements.} N.Cramp\'e and L.Poulain d'Andecy are partially supported by Agence National de la Recherche Projet AHA ANR-18-CE40-0001.
The research of L.Vinet is supported in part by a Discovery Grant from the Natural Science and Engineering Research Council (NSERC) of Canada.
 L.Poulain d'Andecy warmly thanks the Centre de Recherches Math\'ematiques (CRM) for support during his visit to Montreal in the course of this investigation. 
 N. Cramp\'e thanks the Universit\'e de Reims Champagne-Ardenne for its hospitality.

\appendix

\section{Calculation of  $X_{\tm}$ and $Y_{\tm}$ \label{app:XY}}

In this appendix, a basis of the multiplicity space $M_\tm$ is given in terms of the
highest-weight vectors of weight $(m''_1,m''_2)$ in the tensor product $[m_1,m_2]\otimes [m'_1,m'_2]$.
The general idea is borrowed from \cite{PST} but, let us recall that the basis constructed here is different and 
has the advantage of being better adapted to the symmetries. Some details on the calculation 
of the symmetric forms (\ref{eq:XY}) of $X_{\tm}$ and $Y_{\tm}$ are also provided. 

\paragraph{Extremal projector.} As in \cite{PST}, the extremal projector $P_{(m''_1,m''_2)}$ of $\mathfrak{su}_3$ from \cite{AST} is used to build highest-weight vectors.
Its definition, adapted to our setting, and its necessary properties are recalled. 

A basis of the irreducible representation $[m''_1,m''_2]$ consists of weight vectors, and up to a scalar, there is a unique vector of weight $(m''_1,m''_2)$ (the highest-weight vector). 
A projector $P_{(m''_1,m''_2)}$ acting on $[m''_1,m''_2]$ is defined by:
\[P_{(m''_1,m''_2)}(v)=\left\{\begin{array}{ll} v & \text{if $v$ is the highest-weight vector of $[m''_1,m''_2]$,}\\[0.5em]
 0 & \text{if $v$ is a weight vector of weight different from $(m''_1,m''_2)$.}\end{array}\right.\]
The projector $P_{(m''_1,m''_2)}$ is also set to 0 on any irreducible finite-dimensional representation of $\mathfrak{su}_3$ different from $[m''_1,m''_2]$. 
Thus, the element $P_{(m''_1,m''_2)}$ is defined as a linear operator on $[m_1,m_2]\otimes [m'_1,m'_2]$ (since this representation decomposes as a direct sum of 
irreducible finite-dimensional representation of $\mathfrak{su}_3$) and projects the whole space $[m_1,m_2]\otimes [m'_1,m'_2]$ onto the multiplicity space  $M_\tm$. 

From this definition, it is also clear that $P_{(m''_1,m''_2)}$ satisfies:
\begin{equation}\label{eq-P}
P_{(m''_1,m''_2)}.(e_{ij}\otimes 1+1\otimes e_{ij})=0\ \ \ \ \ \text{if $i>j$,}
\end{equation}
where the equality is meant as linear operators on $[m_1,m_2]\otimes [m'_1,m'_2]$. Moreover, recall that any element commuting with the diagonal action of $\mathfrak{su}_3$ on $[m_1,m_2]\otimes [m'_1,m'_2]$ 
sends a highest-weight vector on a linear combination of highest-weight vectors with the same weight. Therefore, $P_{(m''_1,m''_2)}$ commutes with any such element
and, in particular, with $X_{\tm}$ and $Y_{\tm}$.

\paragraph{A generating set of $M_{\tm}$.} 
Let $v_0=v_{(m_1,m_2)}\otimes v_{(m'_1,m'_2)}$, where $v_{(m_1,m_2)}$ (respectively, $v_{(m'_1,m'_2)}$) is a highest-weight vector of $[m_1,m_2]$ (respectively, $[m'_1,m'_2]$).
Since $[m_1,m_2]$ and $[m'_1,m'_2]$ are highest-weight representations, any vector of $[m_1,m_2]\otimes [m'_1,m'_2]$ is a linear combination of vectors of the form:
\begin{equation}\label{vectors}
(e_{31}^{a'}e_{21}^{b'}e_{32}^{c'}\otimes e_{31}^ae_{21}^be_{32}^c)\cdot v_0\,.
\end{equation}
Recall that we are looking for the image of the projector $P_{(m''_1,m''_2)}$ which satisfies (\ref{eq-P}). Therefore it is enough to consider the set of vectors:
\[P_{(m''_1,m''_2)}\cdot\left(1\otimes e_{31}^ae_{21}^be_{32}^c\right)\cdot v_0\,,\ \ \ \ \ \ \ \ a,b,c\in\mathbb{Z}_{\geq 0}\ .\]
An easy calculation shows that $e_{31}e_{21}^b=\frac{1}{b+1}(e_{32}e_{21}^{b+1}-e_{21}^{b+1}e_{32})$ and, by induction, any of the above vectors can be rewritten as a linear combination of the following ones:
\[P_{(m''_1,m''_2)}\cdot\left(1\otimes e_{32}^ae_{21}^be_{32}^c\right)\cdot v_0\,,\ \ \ \ \ \ \ \ a,b,c\in\mathbb{Z}_{\geq 0}\ .\]
Since $P_{(m''_1,m''_2)}$ sends to $0$ any vector of weight different from $(m''_1,m''_2)$, restrict $a,b,c$ such that the vector, before the application of $P_{(m''_1,m''_2)}$, is of weight $(m''_1,m''_2)$. 
This implies that $b=n-1$ and that $a+c=\ell-1$. So, with a choice of normalisation that will be convenient later, we conclude that the space $M_{\tm}$ is spanned by the following vectors:
\begin{equation}\label{wm}
\omega_m:=P_{(m''_1,m''_2)}\cdot \left(1\otimes \frac{e_{32}^m\,e_{21}^{n-1}\,e_{32}^{\ell-1-m}}{m!\,(n-1)!\,(\ell-1-m)!}\right)\cdot v_0\,,\ \ \ \ \ \text{where $m\in\{0,\dots,\ell-1\}$.}
\end{equation}
We can also start with the opposite order on the generators of $\mathfrak{su}_3$ in (\ref{vectors}) and use a similar reasoning exchanging the role of $e_{21}$ and $e_{32}$. For some reasons, it will prove convenient to put this other 
spanning set of $M_{\tm}$ in the following form:
\begin{equation}\label{wpm}\omega'_m:=P_{(m''_1,m''_2)}\cdot \left(\frac{e_{21}^m\,e_{32}^{\ell-1}\,e_{21}^{n-1-m}}{m!\,(\ell-1)!\,(n-1-m)!}\otimes 1\right)\cdot v_0\,,\ \ \ \ \ \text{where $m\in\{0,\dots,n-1\}$.}
\end{equation}
We will use the conventions $\omega_m:=0$ if $m\notin\{0,\dots,\ell-1\}$ and $\omega'_m:=0$ if $m\notin\{0,\dots,n-1\}$.

\paragraph{A basis of $M_{\tm}$.} Let $v\in [m_1,m_2]\otimes [m'_1,m'_2]$ such that $(1\otimes e_{23})\cdot v=0$ and $v$ is an eigenvector of $1\otimes h_{2}$ 
with eigenvalues $N\in\mathbb{Z}_{\geq 0}$. Then $v$ is the highest-weight vector of a copy of the irreducible representation of dimension $N+1$ of the subalgebra generated by $1\otimes e_{32},1\otimes h_{2},1\otimes e_{23}$ isomorphic to $\mathfrak{su}_2$. In particular, we have that $(1\otimes e_{32})^{N+1}\cdot v=0$. We will also use this fact for the other $\mathfrak{su}_2$-subalgebras. We deduce that:
\begin{itemize}
\item $(1\otimes e_{32}^{\ell-1-m})\cdot v_0=0$ unless $\ell-1-m<m'_2$, that is, unless $m>\ell-1-m'_2$.
\item $(e_{32}^{m}\otimes 1)\cdot v_0=0$ unless $m<m_2$.
\item $(1\otimes e_{21}^{n-1}\,e_{32}^{\ell-1-m})\cdot v_0=0$ unless $n-1<m'_1+\ell-1-m$, that is, unless $m<m'_1+\ell-n$.
\item $(e_{21}^{n-1}\,e_{32}^{m}\otimes 1)\cdot v_0=0$ unless $n-1<m_1+m$, that is, unless $m>n-m_1-1$.
\end{itemize}
Using the defining formula (\ref{wm}) for the vectors $\omega_m$ together with the property (\ref{eq-P}) for the projector $P_{(m''_1,m''_2)}$, we deduce that:
\[\omega_m=0\ \ \ \ \text{unless $\text{max}\{\ell-m'_2\,,\ n-m_1\,,\ 0\}\leq m< \text{min}\{\ell\,,\ m_2\,,\ m'_1+\ell-n\}$.}\]
So the number of non-zero vectors $\omega_m$ is at most the difference between the min and the max. If $\ell\leq n$, this number is the dimension of $M_{\tm}$. 
We conclude that if $\ell\leq n$, the non-zero vectors $\omega_m$ form a basis of $M_{\tm}$. 

A similar reasoning with the vectors $\omega'_m$ if $n\leq \ell$ exchanges the role of the parameters as follows: $\ell\leftrightarrow n$ and $(m_1,m_2)\leftrightarrow (m'_2,m'_1)$.
Therefore, a basis of $M_{\tm}$ is
\begin{eqnarray}
\{\omega_m\}_{\text{max}\{\ell-m'_2\,,\ n-m_1\,,\ 0\}\leq m< \text{min}\{\ell\,,\ m_2\,,\ m'_1+\ell-n\}}\ \ \ \ \ \ \text{(if $\ell\leq n$)\,,}\label{basis1}\\
\{\omega'_m\}_{\text{max}\{\ell-m'_2\,,\ n-m_1\,,\ 0\}\leq m< \text{min}\{n\,,\ m'_1\,,\ m_2+n-\ell\}}\ \ \ \ \ \ \text{(if $n\leq \ell$)\,.} \label{basis11}
\end{eqnarray}

\paragraph{Matrix for $X_{\tm}$.} Two equivalent expressions for $X$ are:
\[X=\frac{1}{2}(T^{(1,1,2)}-T^{(1,2,2)})+\frac{1}{3}(\ell_1-\ell_2)=T^{(1,1,2)}+\frac{1}{6}(3l_1-l_2-l_3)+\frac{3}{4}(k_3-k_2-k_1)\,.\]
Recall that $X$ acts on the space $M_{\tm}$, the resulting operator being denoted $X_{\tm}$. We are going to calculate the action of $X$ on the vectors $\omega_m$ defined in (\ref{wm}). 
In fact, due to the second expression of $X$, it would be enough to calculate the action of $T^{(1,1,2)}$, since the remaining terms are Casimir elements which act on $M_{\tm}$ simply by 
numbers given in Remark \ref{rem-Cas}. The strategy presented below is equally applicable for any elements of the centraliser, so for simplicity, we present it for $X$, 
but we mention that the actual technical computations are simpler to perform on $T^{(1,1,2)}$.

We are going to use the following rules:
\begin{equation}\label{rules}
\begin{cases}
 P_{(m''_1,m''_2)}\,(e_{ij}\otimes 1)=-P_{(m''_1,m''_2)}\, (1\otimes e_{ij}) \quad\qquad\text{if $i>j$,}\\[0.5em]
(e_{ij}\otimes 1)\,(1\otimes\gamma_m)\cdot v_0=0 \qquad\quad\hspace{3.4cm} \text{if $i<j$,}\\[0.5em]
(h_1^ah_2^b\otimes h_1^ch_2^d)\,(1\otimes \gamma_m)\cdot v_0=(m_1-1)^a(m_2-1)^b(m'_1+\ell-2n)^c(m'_2+n-2\ell)^d (1\otimes \gamma_m)\cdot v_0,\qquad
\end{cases}
\end{equation}
where $\gamma_m:=0$ if $m\notin\{0,\dots,\ell-1\}$ and $\gamma_m:=  \frac{e_{32}^m\,e_{21}^{n-1}\,e_{32}^{\ell-1-m}}{m!\,(n-1)!\,(\ell-1-m)!}$, otherwise.
The first rule is (\ref{eq-P}), the second rule follows from the fact that $e_{ij}\otimes 1$ will hit $v_0$ which is the tensor product of two highest-weight vectors, 
and the third rule comes from recalling that $\gamma_m\cdot v_0$ is an eigenvector of the Cartan generators and from calculating the eigenvalues.

The strategy is to start by rewriting $X$ as a linear combination of monomials in the generators of $U(\mathfrak{su}_3)\otimes U(\mathfrak{su}_3)$ organized 
(in each factor) in the following order: $e_{21},e_{31},e_{32},e_{12},e_{13},e_{23},h_1,h_2$. Then we use that $X$ commutes with $P$ and write:
\begin{equation}\label{Xomega}
X\cdot \omega_m=XP_{(m''_1,m''_2)}\,(1 \otimes \gamma_m)\cdot v_0=P_{(m''_1,m''_2)}X\,(1 \otimes \gamma_m)\cdot v_0=
P_{(m''_1,m''_2)}\,\sum_{W\in S} x_W (1\otimes W\gamma_m)\,\cdot v_0\ ,
\end{equation}
where $S=\{1,\, e_{21}e_{12},\, e_{32}e_{23},\, e_{31}e_{13},\, e_{32}e_{21}e_{13},\, e_{31}e_{23}e_{12}\}$ and $x_W$ are numbers.
The last equality has been obtained by using only the rules (\ref{rules}). 
Then, for each $W\in S$, $W\gamma_m$ is put in the same order as above and applied on $v_0$. 
Finally, after a tedious but straightforward calculation, one can check that:
\begin{equation}\label{Xomegam}
X\cdot \omega_m=\alpha_m\,\omega_{m-1}+\beta_m\,\omega_m+\delta_m\,\omega_{m+1}\,,
\end{equation}
with $E_{m,k}=\sum_{i=1}^6(\xi_i-m-\frac{1}{2})^k$,
\[\alpha_m=(m-\xi_1)(m-\xi_2)(m-\xi_3)\ ,\ \ \ \ \ \ \ \ \delta_m=(m+1)(m+1-\xi_4)(m+1-\xi_5)\ ,\]
\[\ \ \ \ \ \ \ \ \beta_m=-\frac{1}{108}\Bigl(\frac{7}{2}E_{m,1}^3-18 E_{m,1}E_{m,2}+18 E_{m,3}\Bigr)-\frac{1}{24}E_{m,1}\Bigl(\Lambda^2+2\Bigr)\ ,\]
and $(\xi_1,\,\xi_2,\,\xi_3,\,\xi_4,\,\xi_5,\,\xi_6)$ and $\Lambda$ defined in Section \ref{subsec-formulas} for the situation $\ell\leq n$. 
In the case $\ell\leq n$, we have proved that a basis is given by the subset of vectors $\omega_m$ with indices $\text{max}\{\ell-m'_2\,,\ n-m_1\,,\ 0\}\leq m< \text{min}\{\ell\,,\ m_2\,,\ m'_1+\ell-n\}$, 
and moreover that all other vectors $\omega_m$ are $0$.  
So it remains only to extract the diagonal block of $X$ corresponding to the restricted set of indices to obtain the matrix for $X_{\tm}$. 
This proves the formula given in Section \ref{subsec-formulas} for $X_{\tm}$ in the case $\ell\leq n$.

In the situation $n\leq \ell$, we have to consider the vectors $\omega'_m$ instead of $\omega_m$. We see, comparing the definitions (\ref{wm}) and (\ref{wpm}), that $\omega'_m$ is obtained from $\omega_m$ by applying an automorphism on what is between the projector and $v_0$. This automorphism is the composition of the exchange of the two factors and the automorphism $\sigma$ defined in (\ref{sigma}) (in fact this gives $\omega'_m$ up to a global sign not depending on $m$). Then the key fact is that $X$ is invariant under the same automorphism as indicated in Section \ref{sec:ptt}. From this it is immediate to deduce that the action of $X$ on $\omega'_m$ is as in (\ref{Xomegam}), namely:
\[X\cdot \omega'_m=\alpha'_m\,\omega'_{m-1}+\beta'_m\,\omega'_m+\delta'_m\,\omega'_{m+1}\,,\]
where $\alpha'_m,\beta'_m,\delta'_m$ are obtained from $\alpha_m,\beta_m,\delta_m$ with the following exchange of the parameters: $\ell\leftrightarrow n$ and $(m_1,m_2)\leftrightarrow (m'_2,m'_1)$. So the action of $X$ on the vectors $\omega'_m$, $m=0,\dots,\ell-1$, is also by a tridiagonal matrix and in this case, similarly as above, it remains only to extract the square submatrix corresponding to the indices $\text{max}\{\ell-m'_2\,,\ n-m_1\,,\ 0\}\leq m< \text{min}\{n\,,\ m'_1\,,\ m_2+n-\ell\}$. 
This concludes the verification of the formulas give in Section \ref{subsec-formulas} for $X_{\tm}$ in all cases. 

\paragraph{Matrix for $Y_{\tm}$.}
The strategy used above for calculating the action of $X$ on the vectors $\omega_m$ and $\omega'_m$ is also applicable for the computation of $Y$
rewritten as follows to simplify the computation:
\[ Y = T^{(1,1,2,2)}-\frac{1}{2}(l_1+l_2+l_3)-\frac{1}{48}(k_3-k_1-k_2)^2-\frac{5}{12}k_1k_2-(k_3-k_1-k_2)\ .\]
The calculations are more complicated but leads similarly to all the formulas given in Section \ref{subsec-formulas}.
 
\section{Parameters $z$ of Proposition \ref{lem1} \label{app:B}}

The parameters $z$ in relations \eqref{eq:Xe} and \eqref{eq:Ye} are given by 
\begin{subequations}
\begin{eqnarray}
  z_2 &=& \frac{\xi_3+\xi_4}{2}-\frac{\xi_1+\xi_2+\xi_5+\xi_6}{4}\,,\quad  z_3 =\frac{\xi_3 -\xi_4}{2}\,,\quad  z_1=z_3^2 -\frac{\Lambda^2}{4}\ , \quad z_4=\frac{1}{2}\ ,\\
  z_0&=&\frac{A_3}{4}  -\frac{z_2}{3}\left(\frac{\Lambda^2}{4}  +z_3^2-\frac{2z_2^2}{9}-A_2-\frac{1}{2} \right)     \, ,\\
  z_5 &=&\frac{z_3(A_2\Lambda+z_3)}{3}  -\frac{z_2A_3}{3} -\frac{A_2(\Lambda-2z_3+4)(\Lambda-2z_3+2)}{24}-\frac{(\Lambda^2-4-4z_3^2)^2}{192}  \\
  && -\frac{A_3(\xi_1+\xi_2+\xi_5+\xi_6)}{8}-\frac{(2A_2-1)(\xi_1+\xi_2+\xi_5+\xi_6)^2+(\xi_1^2+\xi_2^2+\xi_5^2+\xi_6^2)^2}{48} -\frac{\xi_1^4+\xi_2^4+\xi_5^4+\xi_6^4}{12}\,,\nonumber\\
 z_6&=& \frac{(4z_3^2-\Lambda^2)z_2}{6}-\frac{A_3}{2}  \,,\quad z_7=\frac{(\Lambda+1)(2z_3-1)}{4}  -\frac{z_2^2}{3}\ , z_8= \frac{( 3\Lambda-2z_3 )z_2}{6} \,,\quad z_9=-\frac{2z_2}{3}\,,\\
z_{10}&=&\frac{(2z_3-\Lambda)(\Lambda+2-2z_3)}{4}\,,\quad z_{11}=\frac{\Lambda}{2}-z_3\,,\quad z_{12}=z_3-\frac{\Lambda}{2}-1\, 
\end{eqnarray}
\end{subequations}
where $A_2$ and $A_3$ are given by \eqref{eq:AA} with $\eta_i$ replaced by $\xi_i$. 
 
  \section{Root poset of $E_6$}\label{app-diagram}
The vertices of the graph are all the positive roots of $E_6$. The edges with label $i$ indicate the action of the simple reflection $s_i$. The notation $1^{c_1}...6^{c_6}$ for roots is a short-hand notation for $\sum_{i=1,\dots,6}c_i\alpha_i$.

In the top left subset of 15 roots, there are $9$ roots between $\alpha_3$ and $12345$ which contains $3$. 
These nine roots are sent ``as a block'' to the right by $s_6$. Similarly the two blocks of underlined roots are exchanged by $s_3$. 
When the action of a simple reflection $s_i$ on a root is not indicated, it means that it leaves this root invariant
with the exception of the action of $s_i$ on $\alpha_i$, which gives $s_i(\alpha_i)=-\alpha_i$. 

With this diagram, it is straightforward to calculate the action of any word in the generators of $W(E_6)$ on any root.
\vskip .2cm
\begin{center}
\begin{tikzpicture}[scale=0.5]
\node at (0,0) {$\alpha_1$};
\node at (4,0) {$\alpha_2$};
\node at (8,0) {$\alpha_3$};
\node at (12,0) {$\alpha_4$};
\node at (16,0) {$\alpha_5$};

\draw (0,-0.4)--(2,-2.6);\draw (4,-0.4)--(2,-2.6);\node[right] at (0.8,-1.3) {\footnotesize{$2$}};\node[left] at (3.2,-1.3) {\footnotesize{$1$}};
\draw (4,-0.4)--(6,-2.6);\draw (8,-0.4)--(6,-2.6);\node[right] at (4.8,-1.3) {\footnotesize{$3$}};\node[left] at (7.2,-1.3) {\footnotesize{$2$}};
\draw (8,-0.4)--(10,-2.6);\draw (12,-0.4)--(10,-2.6);\node[right] at (8.8,-1.3) {\footnotesize{$4$}};\node[left] at (11.2,-1.3) {\footnotesize{$3$}};
\draw (12,-0.4)--(14,-2.6);\draw (16,-0.4)--(14,-2.6);\node[right] at (12.8,-1.3) {\footnotesize{$5$}};\node[left] at (15.2,-1.3) {\footnotesize{$4$}};

\node at (2,-3-0.1) {$12$};
\node at (6,-3-0.1) {$23$};
\node at (10,-3-0.1) {$34$};
\node at (14,-3-0.1) {$45$};

\draw (2,-3.4-0.1)--(4,-5.6-0.1);\draw (6,-3.4-0.1)--(4,-5.6-0.1);\node[right] at (2.8,-4.3-0.1) {\footnotesize{$3$}};\node[left] at (5.2,-4.3-0.1) {\footnotesize{$1$}};
\draw (6,-3.4-0.1)--(8,-5.6-0.1);\draw (10,-3.4-0.1)--(8,-5.6-0.1);\node[right] at (6.8,-4.3-0.1) {\footnotesize{$4$}};\node[left] at (9.2,-4.3-0.1) {\footnotesize{$2$}};
\draw (10,-3.4-0.1)--(12,-5.6-0.1);\draw (14,-3.4-0.1)--(12,-5.6-0.1);\node[right] at (10.8,-4.3-0.1) {\footnotesize{$5$}};\node[left] at (13.2,-4.3-0.1) {\footnotesize{$3$}};

\node at (4,-6-0.2) {$123$};
\node at (8,-6-0.2) {$234$};
\node at (12,-6-0.2) {$345$};

\draw (4,-6.4-0.2)--(6,-8.6-0.2);\draw (8,-6.4-0.2)--(6,-8.6-0.2);\node[right] at (4.8,-7.3-0.2) {\footnotesize{$4$}};\node[left] at (7.2,-7.3-0.2) {\footnotesize{$1$}};
\draw (8,-6.4-0.2)--(10,-8.6-0.2);\draw (12,-6.4-0.2)--(10,-8.6-0.2);\node[right] at (8.8,-7.3-0.2) {\footnotesize{$5$}};\node[left] at (11.2,-7.3-0.2) {\footnotesize{$2$}};

\node at (6,-9-0.3) {$1234$};
\node at (10,-9-0.3) {$2345$};
\draw (6,-9.4-0.3)--(8,-11.6-0.3);\draw (10,-9.4-0.3)--(8,-11.6-0.3);\node[right] at (6.8,-10.3-0.3) {\footnotesize{$5$}};\node[left] at (9.2,-10.3-0.3) {\footnotesize{$1$}};
\node at (8,-12-0.4) {$12345$};

\node at (25,0) {$\alpha_6$};
\draw (25,-0.4)--(25,-2.6);\node[right] at (25,-1.5) {\footnotesize{$3$}};
\node at (25,-3-0.1) {$36$};
\draw (25,-3.4)--(23,-5.6);\node[left] at (24.2,-4.3) {\footnotesize{$2$}};
\draw (25,-3.4)--(27,-5.6);\node[right] at (25.8,-4.3) {\footnotesize{$4$}};

\node at (23,-6-0.2) {$236$};
\node at (27,-6-0.2) {$346$};
\draw (23,-6.4-0.2)--(21,-8.6-0.2);\node[left] at (22.2,-7.3-0.2) {\footnotesize{$1$}};
\draw (23,-6.4-0.2)--(25,-8.6-0.2);\draw (27,-6.4-0.2)--(25,-8.6-0.2);\node[right] at (23.8,-7.3-0.2) {\footnotesize{$4$}};\node[left] at (26.2,-7.3-0.2) {\footnotesize{$2$}};
\draw (27,-6.4-0.2)--(29,-8.6-0.2);\node[right] at (27.8,-7.3-0.2) {\footnotesize{$5$}};

\node at (21,-9-0.3) {$1236$};
\node at (25,-9-0.3) {$\underline{2346}$};
\node at (29,-9-0.3) {$3456$};

\draw (21,-9.4-0.3)--(23,-11.6-0.3);\draw (25,-9.4-0.3)--(23,-11.6-0.3);\node[right] at (21.8,-10.3-0.3) {\footnotesize{$4$}};\node[left] at (24.2,-10.3-0.3) {\footnotesize{$1$}};
\draw (25,-9.4-0.3)--(27,-11.6-0.3);\draw (29,-9.4-0.3)--(27,-11.6-0.3);\node[right] at (25.8,-10.3-0.3) {\footnotesize{$5$}};\node[left] at (28.2,-10.3-0.3) {\footnotesize{$2$}};

\node at (23,-12-0.4) {$\underline{12346}$};
\node at (27,-12-0.4) {$\underline{23456}$};
\draw (23,-12.4-0.4)--(25,-14.6-0.4);\draw (27,-12.4-0.4)--(25,-14.6-0.4);\node[right] at (23.8,-13.3-0.4) {\footnotesize{$5$}};\node[left] at (26.2,-13.3-0.4) {\footnotesize{$1$}};
\node at (25,-15-0.5) {$\underline{123456}$};

\node at (15,-12-0.4) {$\underline{23^246}$};
\draw (15,-12.4-0.4)--(13,-14.6-0.4);\node[left] at (14.2,-13.3-0.4) {\footnotesize{$1$}};
\draw (15,-12.4-0.4)--(17,-14.6-0.4);\node[right] at (15.8,-13.3-0.4) {\footnotesize{$5$}};

\node at (13,-15-0.2-0.4) {$\underline{123^246}$};
\node at (17,-15-0.2-0.4) {$\underline{23^2456}$};
\draw (13,-15.4-0.2-0.4)--(11,-17.6-0.2-0.4);\node[left] at (12.2,-16.3-0.2-0.4) {\footnotesize{$2$}};
\draw (13,-15.4-0.2-0.4)--(15,-17.6-0.2-0.4);\draw (17,-15.4-0.2-0.4)--(15,-17.6-0.2-0.4);\node[right] at (13.8,-16.3-0.2-0.4) {\footnotesize{$5$}};\node[left] at (16.2,-16.3-0.2-0.4) {\footnotesize{$1$}};
\draw (17,-15.4-0.2-0.4)--(19,-17.6-0.2-0.4);\node[right] at (17.8,-16.3-0.2-0.4) {\footnotesize{$4$}};

\node at (11,-18-0.3-0.4) {$12^23^246$};
\node at (15,-18-0.3-0.4) {$\underline{123^2456}$};
\node at (19,-18-0.3-0.4) {$23^24^256$};

\draw (11,-18.4-0.3-0.4)--(13,-20.6-0.3-0.4);\draw (15,-18.4-0.3-0.4)--(13,-20.6-0.3-0.4);\node[right] at (11.8,-19.3-0.3-0.4) {\footnotesize{$5$}};\node[left] at (14.2,-19.3-0.3-0.4) {\footnotesize{$2$}};
\draw (15,-18.4-0.3-0.4)--(17,-20.6-0.3-0.4);\draw (19,-18.4-0.3-0.4)--(17,-20.6-0.3-0.4);\node[right] at (15.8,-19.3-0.3-0.4) {\footnotesize{$4$}};\node[left] at (18.2,-19.3-0.3-0.4) {\footnotesize{$1$}};

\node at (13,-21-0.4-0.4) {$12^23^2456$};
\node at (17,-21-0.4-0.4) {$123^24^256$};
\draw (13,-21.4-0.4-0.4)--(15,-23.6-0.4-0.4);\draw (17,-21.4-0.4-0.4)--(15,-23.6-0.4-0.4);\node[right] at (13.8,-22.3-0.4-0.4) {\footnotesize{$4$}};\node[left] at (16.2,-22.3-0.4-0.4) {\footnotesize{$2$}};
\node at (15,-24-0.5-0.4) {$12^23^24^256$};
\draw(15,-24.4-0.5-0.4)--(15,-26.6-0.5-0.4);\node[right] at (15,-25.5-0.5-0.4) {\footnotesize{$3$}};
\node at (15,-27-0.5-0.5) {$12^23^34^256$};
\draw(15,-27.4-0.5-0.5)--(15,-29.6-0.5-0.5);\node[right] at (15,-28.5-0.5-0.4) {\footnotesize{$6$}};
\node at (15,-30-0.5-0.6) {$12^23^34^256^2=\Theta$};

\draw(15,-7.4)--(18,-8.4);\node[above] at (16.5,-7.9){\footnotesize{$6$}};

\draw(18.5,-14.4)--(21,-13.4);\node[above] at (19.7,-13.9){\footnotesize{$3$}};

\end{tikzpicture}
\end{center}

\end{document}